\documentclass[brochure,12pt]{bourbaki}

\usepackage[T1]{fontenc}
\usepackage{lmodern,amssymb,bm}
\usepackage[all]{xy}
\usepackage{enumerate}

\usepackage[french]{babel}
\usepackage{combelow}

\usepackage{color}

\usepackage[utf8]{inputenc}

\usepackage{wasysym}
\usepackage{tikz}
\usepackage{footmisc}
\usepackage{cleveref}

\usepackage[
backend=biber,
style=authoryear, 
citestyle=authoryear-comp,
maxnames=5
]{biblatex}
\usepackage{csquotes}

\DeclareDelimFormat{nameyeardelim}{\addcomma\space}

\DeclareNameAlias{sortname}{given-family}

\AtEveryBibitem{%
  \clearfield{issn} % Remove issn
  \clearfield{doi} % Remove doi

  \ifentrytype{online}{}{% Remove url except for @online
    \clearfield{url}
  }
}

\renewbibmacro{in:}{%
    \ifentrytype{article}{}{\printtext{\bibstring{in}\intitlepunct}}}

\DeclareFieldFormat[article,periodical,inreference]{number}{\mkbibparens{#1}}
\DeclareFieldFormat[article,periodical,inreference]{volume}{\mkbibbold{#1}}
\renewbibmacro*{volume+number+eid}{%
    \printfield{volume}%
    \setunit*{\addthinspace}% NEW
    \printfield{number}%
    \setunit{\addcomma\space}%
    \printfield{eid}}

\DeclareFieldFormat[article,inbook,incollection]{title}{\enquote{#1}\addcomma} 

\date{Avril 2021}
\bbkannee{72\textsuperscript{e} ann\'ee, 2019--2021}  
\bbknumero{1176}

\title{Progrès récents sur la conjecture de Zagier \\ et le programme de Goncharov}
\subtitle{d'après Goncharov, Rudenko, Gangl, \ldots}

\author{Cl\'{e}ment Dupont}
\address{Institut Montpelli\'{e}rain Alexander Grothendieck,\\ Universit\'{e} de Montpellier, CNRS, Montpellier, France}
\email{clement.dupont@umontpellier.fr}

\newcommand{\ZZ}{\mathbb{Z}}
\newcommand{\QQ}{\mathbb{Q}}
\newcommand{\RR}{\mathbb{R}}
\newcommand{\CC}{\mathbb{C}}

\renewcommand{\to}{\rightarrow}
\newcommand{\To}{\longrightarrow}

\newcommand{\Li}{\operatorname{Li}}
\newcommand{\Cor}{\operatorname{Cor}}

\newcommand{\Ext}{\operatorname{Ext}}
\newcommand{\MT}{\mathrm{MT}}
\newcommand{\DM}{\mathrm{DM}}

\newcommand{\MHS}{\mathrm{MHS}}

\newcommand{\heptagon}{
\begin{tikzpicture}
		\newdimen\R
		\R=4pt
		\draw[xshift=2.5\R] (0:\R) \foreach \x in {51.4286,102.8571,...,359} {
                -- (\x:\R)
            }-- cycle (90:\R)  ;
		\end{tikzpicture}}

\begin{document}

\usetikzlibrary{math}

\maketitle

\section*{Introduction}

	La fonction zêta d'un corps de nombres~$F$ (extension finie du corps des rationnels) a été définie par Dedekind sous la forme
	$$\zeta_F(s) = \sum_{\mathfrak{a}} N(\mathfrak{a})^{-s} \hspace{1cm} (s\in\CC, \,\operatorname{Re}(s)>1)\ ,$$
	où la somme porte sur les idéaux non nuls de l'anneau des entiers $\mathcal{O}_F$, et $N(\mathfrak{a})=|\mathcal{O}_F/\mathfrak{a}|$ est la norme. C'est un invariant arithmétique fondamental d'un corps de nombres, qui généralise la fonction zêta de Riemann $\zeta_\QQ(s)=\zeta(s)$. \textcite{zagierhyperbolic, zagierconjecture} a conjecturé que la valeur spéciale $\zeta_F(n)$, pour un entier $n\geq 2$, s'exprime de manière précise en termes d'évaluations en des éléments de~$F$ du $n$-ième polylogarithme classique
	$$\Li_n(z) = \sum_{k\geq 1}\frac{z^k}{k^n} \hspace{1cm} (z\in\CC, \,|z|<1)\ .$$
	Cette conjecture vise à généraliser la formule analytique du nombre de classes, qui exprime le résidu de la fonction zêta en $s=1$ en termes de logarithmes d'unités.
	
	Les travaux fondateurs de \textcite{borelstable, borelzeta} sur la cohomologie du groupe linéaire permettent de relier $\zeta_F(n)$ au covolume d'un groupe de $K$-théorie de $F$ à l'intérieur d'un espace euclidien via le \emph{régulateur}
	$$K_{2n-1}(F) \To \mathbb{R}^{d_n}\ ,$$
	pour un entier $d_n$ qui dépend de $F$, défini plus bas par \eqref{eq: def dn}. La conjecture de Zagier devient alors un énoncé en $K$-théorie : il s'agit de produire un cocycle explicite, exprimé grace à la fonction $\Li_n$, pour la classe de cohomologie du groupe linéaire qui donne naissance au régulateur. Concrètement, cela consiste à découvrir et organiser les équations fonctionnelles des polylogarithmes, ce qui revient souvent à des considérations non triviales de géométrie projective élémentaire. Après la preuve par \textcite{zagierhyperbolic} du cas $n=2$ de la conjecture\footnote{\textcite{zagierhyperbolic} a découvert et prouvé une version faible du cas $n=2$ de la conjecture par des méthodes de géométrie hyperbolique. La preuve de la version forte, qu'on peut trouver chez \textcite{goncharovgeometry}, est un assemblage de résultats dus à \textcite{blochirvine, blochkyoto}, \textcite{dupontsah}, \textcite{dupontdilog}, et \textcite{suslinbloch}. \label{footnote1}}, cette stratégie a permis à \textcite{goncharovbulletin, goncharovseattle, goncharovgeometry, goncharovICM} de prouver le cas $n=3$. On renvoie aux exposés à ce séminaire d'\textcite{oesterle} et de \textcite{cathelineau} pour des comptes rendus de ces résultats. Des obstructions nouvelles apparaissent pour $n=4$, ce qui explique qu'il ait fallu attendre une vingtaine d'années de plus pour voir une preuve de ce cas de la conjecture de Zagier, par \textcite{goncharovrudenko}. Le premier objectif de ce texte est de décrire les enjeux et les grandes lignes de cette preuve.
	
	Afin de comprendre la nécessité d'ingrédients supplémentaires dans le cas $n\geq 4$, un point de vue motivique sur la conjecture de Zagier est utile. La philosophie des motifs, due à Grothendieck, est celle d'une catégorie qui serait le réceptacle d'une théorie de cohomologie universelle pour les variétés algébriques (concentrant en un objet les groupes de cohomologie singulière, de de Rham, $\ell$-adique, etc.) et où les morphismes seraient liés aux cycles algébriques. Dans la vision de \textcite{beilinsonhigher} et de \textcite{delignelettresoule}, le régulateur devrait être induit par la réalisation de Hodge au niveau des groupes d'extensions via la formule
	\begin{equation}\label{eq: ext K theorie intro}
	\Ext^1_{\MT(F)}(\QQ(-n),\QQ(0)) \simeq K_{2n-1}(F)_\QQ\ .
	\end{equation}
	On a noté $\MT(F)$ la catégorie tannakienne des motifs de Tate mixtes
        sur $F$, qui existe inconditionnellement grâce aux résultats de \textcite{borelstable,borelzeta}
        et aux travaux de \textcite{voevodskytriangulated} et
        \textcite{levinetate}, et pour laquelle la formule ci-dessus est
        vérifiée. Dans ce cadre, la conjecture de Zagier affirme que les
        extensions~\eqref{eq: ext K theorie intro} proviennent toutes de certains motifs \og polylogarithmiques \fg{}, dont les $\Li_n(z)$ sont des périodes. 
	
	Pour $n\leq 3$ la situation est particulièrement favorable, puisqu'on s'attend à ce que \emph{tous} les motifs qui sont des extensions itérées de $\QQ(0)$, $\QQ(-1)$, $\QQ(-2)$, $\QQ(-3)$ proviennent de motifs polylogarithmiques. Pour $n\geq 4$ ce n'est plus le cas et de nouveaux motifs doivent être pris en compte, associés aux polylogarithmes \emph{multiples}
	\begin{equation*}\label{eq: def Li multiple intro}
	\Li_{n_1,\ldots,n_r}(z_1,\ldots,z_r) = \sum_{1\leq k_1<\cdots <k_r} \frac{z_1^{k_1}\cdots z_r^{k_r}}{k_1^{n_1}\cdots k_r^{n_r}}\ \cdot
	\end{equation*}	
	Le problème devient beaucoup plus subtil : du point de vue des
        périodes, il faut montrer que \emph{certaines} combinaisons linéaires
        spéciales de polylogarithmes multiples --- celles qui sont des périodes
        d'extensions~\eqref{eq: ext K theorie intro} ---  s'expriment en termes de polylogarithmes classiques. Dans cette direction, \textcite{ganglfour} a apporté une contribution cruciale en démontrant une équation fonctionnelle, prévue par \textcite{goncharovseattle, goncharovgeometry}, reliant $\Li_{3,1}$ et $\Li_4$. Celle-ci a été redécouverte par \textcite{goncharovrudenko} qui l'incorporent à un contexte général, inspiré par les travaux de \textcite{fockgoncharov} sur le dilogarithme, où la structure des polylogarithmes multiples est organisée par la combinatoire des structures amassées. La structure opéradique des dissections de polygones joue notamment un rôle central, en lien avec le formalisme des corrélateurs, qui sont des variantes des polylogarithmes multiples introduites par \textcite{goncharovcorrelators}.
        
       Au vu de la discussion précédente, les résultats dont il sera question dans ce texte s'inscrivent dans une stratégie vers la conjecture de Zagier découpée en deux étapes distinctes et relativement indépendantes :
       \begin{enumerate}[1)]
       \item La première étape consiste à démontrer une conjecture de Zagier \og faible \fg{}\footnote{Ce terme est parfois utilisé pour qualifier les résultats de \textcite{beilinsondeligne} et de \textcite{dejeu}, qui sont de nature différente.} qui fait intervenir les polylogarithmes \emph{multiples}. De manière imprécise, il s'agit de donner une \og formule \fg{} pour le régulateur de Borel en termes de polylogarithmes multiples. \textcite{goncharovarakelov} a apporté une contribution importante dans cette direction en décrivant le régulateur de Borel via une fonction appelée polylogarithme grassmannien (univalué). \textcite{charltonganglradchenkograssmannian} donnent une formule pour une variante multivaluée de cette fonction, introduite par \textcite{goncharovsimple}, en termes de polylogarithmes multiples. Ces résultats suggèrent qu'on dispose aujourd'hui de tous les ingrédients nécessaires à la réalisation de cette première étape\footnote{Les résultats de \textcite{gerdes} vont dans la même direction puisqu'ils donnent une description d'une partie de la $K$-théorie rationnelle des corps en termes de géométrie projective \emph{linéaire} et donnent donc du poids à la conjecture selon laquelle les extensions \eqref{eq: ext K theorie intro} proviennent toutes de motifs polylogarithmiques multiples.}.
       \item La deuxième étape consiste à passer des polylogarithmes multiples aux polylogarithmes classiques. Cette étape se passe au niveau de versions motiviques de ces fonctions, où une structure (comultiplicative) supplémentaire permet de distinguer les périodes d'extensions \eqref{eq: ext K theorie intro}. Le caractère inexplicite du formalisme motivique amène en fait à travailler avec des versions \og symboliques \fg{} des polylogarithmes multiples, ce qui demande une compréhension fine des relations fonctionnelles qu'ils satisfont.
       \end{enumerate}
	
	\subsection*{Le programme de Goncharov}
	
		Un des enjeux de la conjecture de Zagier est de décrire explicitement la~$K$\nobreakdash-théorie des corps de nombres. Les travaux de \textcite{matsumoto} sur le $K_2$ prolongés par \textcite{milnorktheory}, ainsi que ceux de \textcite{suslinbloch} sur le $K_3$, laissent à penser que le cas des corps de nombres n'est pas spécial et qu'il existe des descriptions explicites uniformes des groupes de $K$-théorie de tous les corps. Le point de vue motivique permet d'éclairer cette question : en admettant des conjectures générales on dispose pour tout corps $F$ d'une catégorie $\MT(F)$ où les groupes d'extensions sont reliés à la~$K$-théorie par une formule qui raffine~\eqref{eq: ext K theorie intro}. Le formalisme tannakien fournit donc des \emph{complexes motiviques} qui calculent les groupes de $K$-théorie de $F$ (ou plus précisément les gradués pour la~$\gamma$-filtration définie par \textcite{souleoperations}). 
		
		Une idée centrale de \textcite{goncharovseattle, goncharovgeometry} est de s'inspirer de ces complexes motiviques pour donner des descriptions inconditionnelles de la $K$-théorie des corps qui soient aussi explicites et \og petites \fg{} que possible. Du point de vue motivique, l'enjeu est de donner une définition inconditionnelle de la catégorie $\MT(F)$ \og par générateurs et relations \fg{},  le rôle des générateurs étant joué par les polylogarithmes multiples (motiviques) ou par certaines variantes. On rassemble ici sous le terme \emph{programme de Goncharov} l'ensemble cohérent de constructions et de conjectures développées par Goncharov dans cette direction durant les trente dernières années. Le second objectif de ce texte est d'en présenter certaines idées directrices (et notamment la \emph{conjecture de liberté} et la \emph{conjecture de profondeur}) et de placer dans ce cadre les résultats de \textcite{ganglfour} et de \textcite{goncharovrudenko}, mais aussi des contributions de \textcite{charltonganglradchenkograssmannian} et de \textcite{rudenkodepth}.

	\paragraph*{Organisation de ce texte}
	
		Au \S 1 on présente la conjecture de Zagier et ses liens avec la $K$-théorie. On y aborde le problème des relations polylogarithmiques et on définit les \emph{complexes polylogarithmiques} de \textcite{goncharovseattle, goncharovgeometry}. Le \S 2 introduit les catégories de motifs de Tate mixtes et leur formalisme tannakien. Dans ce cadre, on étudie dans les \S\S 3-4 les versions motiviques des polylogarithmes (multiples), qui mènent à une interprétation motivique de la conjecture de Zagier. On énonce et analyse la \emph{conjecture de liberté} et la \emph{conjecture de profondeur} de Goncharov, qui sont des énoncés généraux sur la structure de la catégorie des motifs de Tate mixtes sur un corps. 
		
		Après avoir introduit les corrélateurs, qui sont des variantes des polylogarithmes multiples, on étudie au \S 5 les résultats de \textcite{ganglfour} et de \textcite{goncharovrudenko} sur la structure fine des motifs de Tate mixtes en poids $\leq 4$, qui mènent à la preuve de la conjecture de Zagier pour $n=4$. Un rôle important est joué par certaines familles de relations entre polylogarithmes multiples qui trouvent leur origine dans l'étude des structures amassées, point de vue qu'on développe au \S 6. On conclut en mentionnant les travaux récents de \textcite{rudenkodepth} sur la conjecture de profondeur.
		
		Il est ici question d'un sujet aux multiples ramifications
                dont certains des aspects importants n'apparaissent pas dans
                ce texte, par manque de place. C'est notamment le cas des
                liens avec la géométrie hyperbolique (voir par exemple \textcite{goncharovvolumes}, \textcite{browndedekind}, \textcite{rudenkodepth}) et des aspects \og non génériques \fg{} reliés à la $K$-théorie de variétés générales et aux classes caractéristiques (voir par exemple \textcite{goncharovexplicit}, \textcite{goncharovrudenko}).
	
	\paragraph*{Notations et conventions}
	Pour des nombres complexes non nuls~$a,b$ on note $a\sim_{\QQ^\times}b$ pour signifier que $a\in \QQ^\times b$. Pour un groupe abélien~$A$ on note $A_\QQ=A\otimes_\ZZ\QQ$. On note~$\QQ[X]$ le $\QQ$-espace vectoriel librement engendré par les éléments d'un ensemble~$X$. Pour une application $f$ de $X$ vers un $\QQ$-espace vectoriel $V$, on note $f:\QQ[X]\to V$ l'application linéaire induite. On note~$V_n$ la composante homogène de degré~$n$ d'un espace vectoriel $\ZZ$-gradué $V$. Tous les espaces vectoriels, algèbres, etc. sont implicitement définis sur $\mathbb{Q}$.
	
	\paragraph*{Remerciements}
	Je tiens à remercier Herbert Gangl, Alexander Goncharov, et Daniil Rudenko pour leur aide précieuse dans la préparation de ce texte et la patience avec laquelle ils m'ont expliqué leurs travaux. Un grand merci aussi à Francis Brown, José Ignacio Burgos Gil, Damien Calaque, Rob de Jeu, Javier Fres\'{a}n, Richard Hain, et Don Zagier pour leurs commentaires et suggestions sur une version préliminaire de ce texte.

\section{La conjecture de Zagier}

	\subsection{La formule analytique du nombre de classes}
	
		Soit $F$ un corps de nombres, dont on note $r_1$ le nombre de plongements réels (numérotés $\sigma_1,\ldots,\sigma_{r_1}$) et $r_2$ le nombre de paires de plongements complexes conjugués non réels (numérotés $\sigma_{r_1+1},\ldots,\sigma_{r_1+r_2},\overline{\sigma_{r_1+1}},\ldots,\overline{\sigma_{r_1+r_2}}$), de sorte que le degré de $F$ est~$[F:\QQ]=r_1+2r_2$. Le \emph{régulateur de Dirichlet} est le morphisme de groupes abéliens
		\begin{equation}\label{eq: regulateur dirichlet}
		\rho:\mathcal{O}_F^\times\oplus \ZZ\To \RR^{r_1+r_2}
		\end{equation}
		défini par $\rho(x)=(\log|\sigma_1(x)|,\ldots,\log|\sigma_{r_1}(x)|, \log|\sigma_{r_1+1}(x)|^2,\ldots,\log|\sigma_{r_1+r_2}(x)|^2)$
                pour $x\in\mathcal{O}_F^\times$, et $\rho(k)=\frac{1}{r_1+r_2}(k,\ldots,k)$
                pour~\mbox{$k\in\mathbb{Z}$}. 
		Le théorème des unités de Dirichlet affirme que $\rho$ est injectif modulo torsion et que son image est un réseau de $\RR^{r_1+r_2}$ ; le régulateur de $F$, noté $R_F$, est le covolume de ce réseau.
		
		La \emph{formule analytique du nombre de classes} exprime le résidu de la fonction zêta de~$F$ en $s=1$ en termes du régulateur, du nombre~$w_F$ de racines de l'unité dans~$F$, du discriminant~$D_F$, et du nombre de classes~$h_F$ :
		\begin{equation}\label{eq: formule nombre classes}
		\lim_{s\to 1}\,(s-1)\zeta_F(s) = \frac{2^{r_1+r_2}h_F}{w_F}\frac{\pi^{r_2}}{\sqrt{|D_F|}} R_F\ .
		\end{equation}
		Comme son nom l'indique, son intérêt premier est de donner accès au nombre de classes. Cependant, nous nous intéresserons ici à la partie \emph{transcendante} de la formule, que nous écrirons sous la forme
		\begin{equation}\label{eq: formule nombre classes simplifiee}
		\lim_{s\to 1}\,(s-1)\zeta_F(s) \sim_{\QQ^\times} \frac{\pi^{r_2}}{\sqrt{|D_F|}} \det(\log|\sigma_i(\varepsilon	_j)|)_{1\leq i,j\leq r_1+r_2-1}\ ,
		\end{equation}
		où $(\varepsilon_1,\ldots,\varepsilon_{r_1+r_2-1})$ est une base de $\mathcal{O}_F^\times$ modulo torsion.

	\subsection{Polylogarithmes et conjecture de Zagier}
	
		\textcite{zagierconjecture} a conjecturé une généralisation de la formule~\eqref{eq: formule nombre classes simplifiee} à toutes les valeurs spéciales des fonctions zêta de Dedekind, où le logarithme doit être remplacé par les polylogarithmes. Pour un entier $n\geq 1$, rappelons la définition du $n$-ième \emph{polylogarithme classique} :
		\begin{equation}\label{eq: def Li n}
		\Li_n(z) = \sum_{k\geq 1}\frac{z^k}{k^n} \hspace{1cm} (z\in\CC, \,|z|<1)\ .
		\end{equation}
		On renvoie le lecteur à \textcite{lewin} pour une histoire de
                ces fonctions qui remonte à Leibniz. On a
                $\Li_1(z)=-\log(1-z)$ et l'équation différentielle
                $\Li_n'(z)=\frac{1}{z}\Li_{n-1}(z)$, pour~$n\geq 2$,  montre que \eqref{eq: def Li n} se prolonge en une fonction holomorphe \emph{multivaluée} sur le domaine $\mathbb{C}\setminus\{1\}$. Le caractère multivalué signifie que la détermination de $\Li_n(z)$ change quand on fait un prolongement en suivant un lacet autour de $z=1$. Zagier a introduit les versions \emph{univaluées}
		$$\mathrm{P}_n(z) = p_n\left(\sum_{k=0}^{n-1} \frac{2^kB_k}{k!} \log^k|z| \Li_{n-k}(z)\right) = p_n\big(\Li_n(z) - \log|z|\Li_{n-1}(z) + \cdots\big)\ ,$$
		où les $B_k$ sont les nombres de Bernoulli\footnote{C'est-à-dire $B_0=1$, $B_1=-\frac{1}{2}$, $B_2=\frac{1}{6}$, $B_3=0$, $B_4=-\frac{1}{30}$, etc. } et où $p_n:\CC\to \CC/(2\mathrm{i}\pi)^n\RR\simeq \RR$ est l'application \og partie réelle \fg{} si $n$ est impair et \og partie imaginaire \fg{} si $n$ est pair. 
		Ce sont des fonctions sur $\mathbb{P}^1(\CC)\setminus \{0,1,\infty\}$ qui sont analytiques réelles et se prolongent par continuité à $\mathbb{P}^1(\CC)$ si $n\geq 2$, avec $\mathrm{P}_n(0)=\mathrm{P}_n(\infty)=0$. Comme on le verra plus bas, les fonctions $\mathrm{P}_n$ sont naturelles du point de vue motivique \parencite{beilinsondeligne}, ce qui explique qu'elles se comportent mieux que d'autres versions univaluées des polylogarithmes classiques (considérées par exemple par \textcite{ramakrishnananalogs} et par \textcite{zagierblochwignerramakrishnan}).
		
		On a $\mathrm{P}_1(z)=-\log|1-z|$, et $\mathrm{P}_2(z)=\mathrm{Im}(\Li_2(z)+\log|z|\log(1-z))$
		est une fonction spéciale remarquable qui est connue dans la littérature sous le nom de \emph{dilogarithme de Bloch--Wigner} \parencite{blochkyoto, blochirvine}. Elle vérifie une équation fonctionnelle \og à $5$ termes \fg{}, pour $x,y\in \mathbb{P}^1(\CC)\setminus \{0,1,\infty\}$ :
		\begin{equation}\label{eq: five term}
		\mathrm{P}_2(x)-\mathrm{P}_2(y)+\mathrm{P}_2(\textstyle\frac{y}{x}) - \mathrm{P}_2(\textstyle\frac{1-y}{1-x}) + \mathrm{P}_2(\textstyle\frac{x(1-y)}{y(1-x)}) = 0
		\end{equation}
		qu'on peut prouver par différentiation en notant\footnote{On peut aussi (voir par exemple \textcite{zagierdilogarithm}) la prouver par une méthode de géométrie hyperbolique en remarquant que $\mathrm{P}_2(x)$ est le volume du tétraèdre hyperbolique \og idéal \fg{} dont les sommets sont $\infty$, $0$, $1$, et $x$ dans le bord de l'espace hyperbolique $\mathbb{H}^3$, identifié à $\mathbb{P}^1(\CC)$. L'équation à $5$ termes exprime l'additivité du volume pour les $5$ tétraèdres hyperboliques dont les sommets sont pris parmi $\infty$, $0$, $1$, $x$, $y$. Il s'agit d'une version univaluée d'une équation fonctionnelle pour $\Li_2(z)$ attribuée à Spence et à Abel.} que $d\mathrm{P}_2(z) = -\log|1-z|\,d\operatorname{arg}(z) + \log|z|\,d\operatorname{arg}(1-z)$ et que le membre de gauche de \eqref{eq: five term} s'annule pour $x=y$.

                Soient $F$ un corps de nombres et $r_1,r_2$ comme dans le paragraphe précédent. Pour un entier $n\geq 2$ posons 
               	\begin{equation}\label{eq: def dn}
               	d_n= \begin{cases} r_1+r_2 & \hbox{si } n \hbox{ est impair} \\ r_2 & \hbox{si } n \hbox{ est pair} \end{cases}
               	\end{equation} 
               	et fixons l'identification, induite par la partie réelle ou la partie imaginaire,
		\begin{equation}\label{eq: identification Rdn}
		\left(\bigoplus_{\sigma:F\to \CC}\CC/(2\mathrm{i}\pi)^n\RR\right)^+ \simeq \RR^{d_n}\ ,
		\end{equation}
		où le symbole $+$ désigne l'espace des invariants pour la conjugaison complexe, qui agit à la fois sur $\mathbb{C}/(2 \mathrm{i}\pi)^n\RR$ et sur l'ensemble des plongements de $F$ dans $\CC$. Si $\mathrm{P}_n$ est vu comme étant à valeurs dans $\CC/(2 \mathrm{i}\pi)^n\RR$ il vérifie $\mathrm{P}_n(\overline{z})=\overline{\mathrm{P}_n(z)}$ et on obtient grâce à~\eqref{eq: identification Rdn} un morphisme noté
		\begin{equation}\label{eq: Ln cal corps de nombres}
		\mathrm{P}_n^F : \QQ[F^\times] \to \RR^{d_n}\ .
		\end{equation}
		Il se calcule, pour $x\in F^\times$, par
		$$
        \mathrm{P}_n^F(x)= \begin{cases} (\mathrm{P}_n(\sigma_1(x)),\ldots,\mathrm{P}_n(\sigma_{r_1+r_2}(x))) & \hbox{si } n \hbox{ est impair} \\  (\mathrm{P}_n(\sigma_{r_1+1}(x)),\ldots,\mathrm{P}_n(\sigma_{r_1+r_2}(x))) & \hbox{si } n \hbox{ est pair} \end{cases}
        $$
		et est étendu par linéarité à $\QQ[F^\times]$.

		\begin{conj}[Zagier]\label{conj: zagier}
        Soient $F$ un corps de nombres et $n$ un entier $\geq 2$. Il existe des éléments $\xi_1,\ldots,\xi_{d_n}\in \QQ[F^\times]$ tels qu'on ait 
		\begin{equation}\label{eq: zagier conj}
		\zeta_F(n) \sim_{\mathbb{Q}^\times}\frac{\pi^{n([F:\QQ]-d_n)}}{\sqrt{|D_F|}} \det\left(\mathrm{P}_n^F(\xi_j)\right)_{1\leq j\leq d_n}\ .
		\end{equation}
		\end{conj}
		
		Dans le cas où $F$ est totalement réel et $n$ est pair, cette conjecture est un résultat de \textcite{klingen} et \textcite{siegel} (on a alors $d_n=0$ et le déterminant dans~\eqref{eq: zagier conj} vaut $1$) qui généralise la formule d'Euler : $\zeta(n)\sim_{\QQ^\times}\pi^n$ pour $n$ pair. Elle est aussi connue pour tout $n$ dans le cas où $F$ est un corps cyclotomique : cela découle en effet de la décomposition de la fonction zêta de $F$ en produit de fonctions $L$ de Dirichlet, et de ce que la valeur en $n$ d'une telle fonction $L$ s'exprime comme combinaison linéaire d'évaluations de $\Li_n$ aux racines de l'unité. Plus généralement, le théorème de Kronecker--Weber permet d'obtenir, pour un corps de nombres abélien $F$, une expression de $\zeta_F(n)$ en termes d'évaluations de $\mathrm{Li}_n$ aux racines de l'unité --- qui ne sont cependant pas des éléments de $F$ en général. Notons enfin que la conjecture \ref{conj: zagier} se prête à (et est justifiée par) des vérifications numériques, comme l'expliquent \textcite{zagierconjecture} et \textcite{zagiergangl}.	
		
		La conjecture \ref{conj: zagier} a été prouvée pour $n=2$ par Zagier (voir la note de bas de page \footref{footnote1}) et pour $n=3$ par \textcite{goncharovbulletin, goncharovseattle, goncharovgeometry}. Le théorème suivant a été prouvé récemment par \textcite{goncharovrudenko}.
		
		\begin{theo}
		La conjecture de Zagier est vraie pour $n=4$.
		\end{theo}

		\begin{rema}\label{rema: equation fonctionnelle}
		En utilisant l'équation fonctionnelle reliant $\zeta_F(s)$ à $\zeta_F(1-s)$, on peut exprimer~\eqref{eq: zagier conj} sous la forme $\zeta_F^*(1-n)\sim_{\QQ^\times}\pi^{(1-n)d_n}\det\left(\mathrm{P}_n(\xi_j)\right)_{1\leq j\leq d_n}$, où $\zeta_F^*(1-n)$ est le premier coefficient non nul dans le développement limité de $\zeta_F(s)$ en $s=1-n$ (qui s'avère être le coefficient de degré $d_n$).
		\end{rema}
		
	\subsection{Polylogarithmes multiples}
	
		Pour des entiers $n_1,\ldots,n_r\geq 1$ on peut généraliser~\eqref{eq: def Li n} et définir en suivant \textcite{goncharovICM} le \emph{polylogarithme multiple} (à plusieurs variables)
		\begin{equation}\label{eq: def Li multiple}
		\Li_{n_1,\ldots,n_r}(z_1,\ldots,z_r) = \sum_{1\leq k_1<\cdots <k_r} \frac{z_1^{k_1}\cdots z_r^{k_r}}{k_1^{n_1}\cdots k_r^{n_r}}\ ,
		\end{equation}
		qui est une fonction holomorphe au voisinage de $z_i=0$ et qu'on peut prolonger analytiquement. L'étude de telles fonctions, aussi appelées \emph{hyperlogarithmes}, remonte au moins à \textcite{kummeriterated} et à \textcite{poincaregroupes}.
		Elles apparaissent naturellement lorsqu'on multiplie des polylogarithmes classiques, par exemple
		\begin{equation*}\label{eq: relation quasi melange}
		\Li_{n_1}(z_1)\Li_{n_2}(z_2) = \Li_{n_1,n_2}(z_1,z_2)+\Li_{n_2,n_1}(z_2,z_1) + \Li_{n_1+n_2}(z_1z_2)\ ,
		\end{equation*}
		qui s'obtient facilement en découpant le domaine de sommation double. On appelle \emph{poids} la somme $n_1+\cdots+n_r$ des indices dans~\eqref{eq: def Li multiple} et \emph{profondeur} le nombre~$r$ d'indices. Un fait général est que les relations fonctionnelles linéaires entre polylogarithmes multiples sont homogènes pour le poids.
		En revanche, pour un poids donné, il existe de nombreuses relations entre polylogarithmes multiples de profondeurs différentes. Par exemple, la relation suivante, valable autour de $x=y=0$ et qui se montre par différentiation, permet d'exprimer la fonction de deux variables $\Li_{1,1}$ comme combinaison linéaire d'évaluations de la fonction d'une variable~$\Li_2$ :
		$$\Li_{1,1}(x,y) = \Li_2(\textstyle\frac{y(1-x)}{y-1}) - \Li_2(\textstyle\frac{y}{y-1}) - \Li_2(xy)\ .$$
		De même, en poids $3$, les fonctions $\Li_{2,1}$, $\Li_{1,2}$ et $\Li_{1,1,1}$ peuvent toutes s'exprimer comme combinaisons linéaires d'évaluations de $\Li_3$ et de produits d'évaluations de $\Li_1$ et $\Li_2$. 
		
		Ce phénomène s'arrête en poids $4$, où $\Li_1$, $\Li_2$, $\Li_3$, $\Li_4$ ne suffisent pas à exprimer tous les polylogarithmes multiples, comme l'a remarqué \textcite{boehm} dans le contexte du calcul des volumes hyperboliques en dimension 7 (voir aussi \textcite{wechsung}, \textcite{wojtkowiakbasic}, et \textcite{goncharovseattle}). Comme nous le verrons plus bas, c'est ce qui rend le cas $n=4$ de la conjecture de Zagier différent des cas $n=2,3$ et crucial dans la compréhension du cas général.
		
	\subsection{$K$-théorie et régulateurs supérieurs}
	
		La conjecture de Zagier est la partie émergée d'un iceberg qui interroge la description de la $K$-théorie des corps. D'après \textcite{quillenK}, les foncteurs de $K$-théorie associent, pour tout entier $i\geq 0$, un groupe abélien $K_i(R)$ à un anneau (unitaire) $R$, qui est un invariant \og homotopique \fg{} de la catégorie des $R$-modules projectifs. Nous serons surtout intéressés par la $K$-théorie rationnelle $K_i(R)_\QQ$, qui est plus facile à calculer. En effet, le théorème de \textcite{milnormoore} 
		permet d'en donner une description \og homologique \fg{} comme l'espace des éléments primitifs de l'homologie rationnelle du groupe général linéaire stable $\mathrm{GL}(R)$ (l'union des groupes $\mathrm{GL}_N(R)$, où l'inclusion $\mathrm{GL}_N(R)\hookrightarrow \mathrm{GL}_{N+1}(R)$ est donnée par $g\mapsto  \left(\begin{smallmatrix}g & 0\\ 0 &1\end{smallmatrix}\right)$) :
		$$K_i(R)_\QQ \simeq \operatorname{Prim}H_i(\mathrm{GL}(R),\QQ)\ .$$
		
		La $K$-théorie rationnelle d'un corps de nombres $F$ a été complètement calculée par \textcite{borelstable, borelzeta}. En plus des calculs classiques $K_0(F)=\mathbb{Z}$ et $K_1(F)=F^\times$, Borel démontre l'annulation $K_{2i}(F)_\QQ=0$ pour $i\geq 1$.  Concernant le degré impair, il définit pour tout $n\geq 2$ une classe de cohomologie dans $H^{2n-1}(\mathrm{GL}(\CC),\CC/(2\mathrm{i}\pi)^n\RR)$, qui induit un morphisme $K_{2n-1}(\CC)\to \CC/(2\mathrm{i}\pi)^n\RR$. Celui-ci est compatible à la conjugaison complexe et induit par fonctorialité et l'identification~\eqref{eq: identification Rdn} le \emph{régulateur de Borel}\footnote{Dans certaines références, et notamment chez \textcite{burgosregulators}, l'identification \eqref{eq: identification Rdn} est remplacée par l'identification induite par $\CC\ni z\mapsto \mathrm{Re}((2\mathrm{i}\pi)^{1-n}z)$, ce qui  multiplie le régulateur de Borel par $\pm (2\pi)^{1-n}$. Cette normalisation a l'effet de simplifier certaines formules ; par exemple, \eqref{eq: borel theorem} s'exprime sous la forme $\zeta_F^*(1-n)\sim_{\QQ^\times}R_F^{(n)}$ (voir la remarque \ref{rema: equation fonctionnelle}). La normalisation que nous choisissons ici est plus adaptée à la conjecture de Zagier.}
		$$\rho_n:K_{2n-1}(F)\To \RR^{d_n}\ .$$
		Il joue pour les groupes de $K$-théorie supérieurs le même rôle que le régulateur de Dirichlet~\eqref{eq: regulateur dirichlet} pour le groupe $K_1(\mathcal{O}_F)\simeq \mathcal{O}_F^\times$ (on a $K_i(\mathcal{O}_F)_\QQ\simeq K_i(F)_\QQ$ pour $i\geq 2$), comme le montre le théorème suivant.
	
		\begin{theo}[Borel]\label{theo: borel}
		Soit $n\geq 2$. Le régulateur $\rho_n$ est injectif modulo torsion et son image est un réseau de $\RR^{d_n}$ (d'où $\dim_\QQ K_{2n-1}(F)_\QQ=d_n$). Le covolume $R_F^{(n)}$ de ce réseau est relié à la fonction zêta de $F$ par la formule
		\begin{equation}\label{eq: borel theorem}
		\zeta_F(n) \sim_{\QQ^\times} \frac{\pi^{n([F:\QQ]-d_n)}}{\sqrt{|D_F|}} R_F^{(n)}\ .
		\end{equation}
		\end{theo}
		
		Les preuves déjà mentionnées de la conjecture de Zagier pour $n=2,3,4$ s'appuient sur le théorème de Borel et la preuve de la conjecture suivante, que nous préciserons plus bas, dont l'enjeu est de donner une \og formule \fg{} pour le régulateur $\rho_n$.
		
		\begin{conj}\label{conj: zagier k theorie}
		Il existe
		\begin{enumerate}[---]
		\item un sous-espace $\mathcal{R}_n(F)\subset \QQ[F^\times]$ ;
		\item une application linéaire $\varphi_n:K_{2n-1}(F)_\QQ\To \mathcal{B}_n(F) := \QQ[F^\times]/\mathcal{R}_n(F)$ ;
		\end{enumerate}
		tels que $\mathrm{P}_n^F$ passe au quotient par $\mathcal{R}_n(F)$ et qu'on ait le diagramme commutatif suivant :
		\begin{equation}\label{eq: diag comm zagier k theorie}
		\xymatrixcolsep{5pc}\xymatrix{
		K_{2n-1}(F)_\QQ \ar[r]^-{\varphi_n} \ar@/_2pc/[rr]_-{\rho_n} & \mathcal{B}_n(F) \ar[r]^-{\mathrm{P}_n^F} & \RR^{d_n}
		}
		\end{equation}
		\end{conj}
		
		En effet, \eqref{eq: zagier conj} découle alors de~\eqref{eq: borel theorem} et~\eqref{eq: diag comm zagier k theorie} en prenant pour $(\xi_j)_{1\leq j\leq d_n}$ des représentants des images par $\varphi_n$ des éléments d'une base de $K_{2n-1}(F)_\QQ$.
		
		\begin{rema}
		Vu que le régulateur de Borel est d'abord défini sur la $K$-théorie de~$\CC$ puis induit sur celle de $F$ par fonctorialité, il est naturel de sortir la conjecture~\ref{conj: zagier k theorie} du cadre des corps de nombres et de vouloir définir $\mathcal{R}_n(F)$ et $\varphi_n$ pour un corps~$F$ quelconque, de manière fonctorielle en $F$ --- ce que nous ferons dans le prochain paragraphe.
		\end{rema}		
		
		\begin{rema}
		Le nombre rationnel implicite dans~\eqref{eq: borel theorem} devrait essentiellement être, selon une conjecture de \textcite{lichtenbaum}, le quotient des ordres des sous-groupes de torsion de $K_i(\mathcal{O}_F)$ pour $i\in \{2n-2,2n-1\}$. Cela généralise naturellement \eqref{eq: formule nombre classes} puisque les sous-groupes de torsion de $K_0(\mathcal{O}_F)$ et $K_1(\mathcal{O}_F)$ sont respectivement le groupe des classes de~$F$ et le groupe des racines de l'unité dans $F$. Cette conjecture est connue pour~$F$ totalement réel et $n$ pair (on renvoie le lecteur à l'article de survol de \textcite{kahnhandbook} pour plus de détails). Dans l'esprit de la formule analytique du nombre de classes \eqref{eq: formule nombre classes}, on peut alors l'appliquer pour calculer l'ordre de la torsion en $K$-théorie (voir par exemple \textcite{burnsetal}).
		\end{rema}
		
		\begin{rema}\label{rema: regulateurs}
		\textcite{beilinsonhigher} a défini des régulateurs pour toutes les variétés algébriques sur $\QQ$ (la construction de Beilinson dans le cas de la variété $\mathrm{Spec}(F)$ redonne le régulateur de Borel multiplié par $\frac{1}{2}$ d'après \textcite{burgosregulators}) et conjecturé une vaste généralisation du théorème de Borel dans ce cadre. Une version de la conjecture de Zagier pour la valeur en $s=2$ de la fonction $L$ d'une courbe elliptique sur $\QQ$, qui fait intervenir le dilogarithme elliptique de \textcite{blochirvine}, a été démontrée par \textcite{goncharovlevin}. 
		\end{rema}
		
	\subsection{Relations polylogarithmiques}
	
		On appelle $\mathcal{R}_n(F)\subset \QQ[F^\times]$ comme dans
                la conjecture~\ref{conj: zagier k theorie} un espace de \emph{relations polylogarithmiques}. \textcite{zagierconjecture} et \textcite{goncharovseattle, goncharovgeometry} ont défini un candidat pour cet espace qui a l'avantage d'avoir un sens pour un corps $F$ quelconque et que nous décrivons maintenant\footnote{Zagier considère également des variantes spécifiques au cas des corps de nombres. Zagier comme Goncharov considèrent en fait des versions entières de $\mathcal{R}_n(F)$ et $\mathcal{B}_n(F)$, respectivement sous-objets et quotients du groupe abélien libre~$\mathbb{Z}[F^\times]$.}.
	
		\subsubsection{Définition récursive}\label{par: def Bcal}
	
			Pour démontrer des équations fonctionnelles pour la fonction $\mathrm{Li}_n$, l'outil principal est de se ramener à la fonction $\mathrm{Li}_{n-1}$ grâce à l'équation différentielle
			\begin{equation}\label{eq: equa diff Li}
			d\mathrm{Li}_n(z) = \mathrm{Li}_{n-1}(z) \,d\log(z)\ .
			\end{equation}
			Plus précisément, on se ramène à des spécialisations d'équations fonctionnelles en une variable $\sum_ia_i\Li_n(x_i(t))=\mathrm{constante}$, avec $a_i\in \QQ$ et $x_i(t)\in \CC(t)$, qui se démontrent par différentiation par rapport à $t$.
			La définition de
			$\mathcal{R}_n(F)$ s'obtient en abstrayant ce procédé.
	
			On définit récursivement sur $n\geq 2$ un sous-espace $\mathcal{R}_n(F)\subset \QQ[F^\times]$ et le quotient 
			$$\mathcal{B}_n(F)=\QQ[F^\times]/\mathcal{R}_n(F)\ , $$ 
			de manière fonctorielle en $F$ ; on note $[x]_n$ la classe de~$x\in F^\times$ dans $\mathcal{B}_n(F)$. Les $\mathcal{B}_n(F)$ sont parfois appelés \emph{groupes de Bloch supérieurs}\footnote{Cette terminologie est parfois utilisée dans la littérature pour désigner d'autres objets.}. Dans tout ce qui suit il sera pratique de voir $\QQ[F^\times]$ comme le quotient de $\QQ[\mathbb{P}^1(F)]$ par les relations $[\infty]=[0]=0$. 
			On considère 
				\begin{equation}\label{eq: delta deux tilde}
				\widetilde{\delta}_2:\QQ[F^\times]\to
                        \Lambda^2F^\times_\QQ\; , \; [x]\mapsto -(1-x)\wedge x
                        \; , \; [1]\mapsto 0\ .
                        \end{equation} 
                        Pour $n\geq 3$, en supposant le foncteur
                        $\mathcal{B}_{n-1}$
                        défini, on considère
                        \begin{equation}\label{eq: delta tilde}
                        \widetilde{\delta}_n :\QQ[F^\times] \to
                        \mathcal{B}_{n-1}(F)\otimes F^\times_\QQ\; ,\;
                        [x]\mapsto [x]_{n-1}\otimes x\ .
                        \end{equation} 
                        Ces formules sont
                        censées rappeler~\eqref{eq: equa diff Li}. On définit
                        alors $\mathcal{R}_n(F)$ comme le sous-espace de
                        $\QQ[F^\times]$ engendré par les différences
                        $\xi(1)-\xi(0)$, pour $\xi(t)\in \QQ[F(t)^\times]$ dans le noyau du morphisme $\widetilde{\delta}_n$ associé au corps $F(t)$.
			
						On montre que la fonction $\mathrm{P}_n$ passe au quotient par $\mathcal{R}_n(\CC)$ et induit $\mathrm{P}_n:\mathcal{B}_n(\CC)\to \RR$ et donc~$\mathrm{P}_n^F:\mathcal{B}_n(F)\to \RR^{d_n}$ pour $F$ un corps de nombres.
			
			\begin{exem}\label{exem: relations dilog}
			Calculons pour $x\in F^\times\setminus\{1\}$ : 
			$$\widetilde{\delta}_2([x]+[x^{-1}]) = -(1-x)\wedge x - (1-x^{-1})\wedge x^{-1} = \frac{1}{1-x}\wedge x+ \frac{x-1}{x}\wedge x = \frac{-1}{x}\wedge x = 0\ .$$
			Ainsi, l'élément $\xi(t)=[tx]+[(tx)^{-1}]$ est dans le noyau de $\widetilde{\delta}_2:\QQ[F(t)^\times]\to \Lambda^2(F(t)^\times_\QQ)$ et l'élément $\xi(1)-\xi(0)=[x]+[x^{-1}]$ est dans $\mathcal{R}_2(F)$. Ce raisonnement est analogue à la preuve de l'équation fonctionnelle $\mathrm{P}_2(z)+\mathrm{P}_2(z^{-1})=0$ par dérivation et évaluation en $z=0$. De même, on montre par récurrence sur $n\geq 2$ que l'élément $[x]_n+(-1)^n[x^{-1}]_n$ est dans $\mathcal{R}_n(F)$, ce qui est analogue à l'équation fonctionnelle $\mathrm{P}_n(z)+(-1)^n\mathrm{P}_n(z^{-1})=0$.
			\end{exem}
			
			\begin{rema}
			On espère pouvoir donner une description concrète (éventuellement conjecturale) de~$\mathcal{R}_n(F)$ en produisant des familles assez générales de relations polylogarithmiques.
			\end{rema}

		\subsubsection{Relations polylogarithmiques pour $n=2,3$}\label{par: relations polylog 2 3}
			
			Un rapide calcul montre que la relation suivante, qui rappelle la relation à $5$ termes \eqref{eq: five term}, est vérifiée dans $\mathcal{B}_2(F)$ pour tous les $x,y\in F\setminus \{0,1\}$ :
			\begin{equation}\label{eq: five term abstract}
			[x]_2 - [y]_2 +\left[\frac{y}{x}\right]_2 - \left[\frac{1-y}{1-x}\right]_2 + \left[\frac{x(1-y)}{y(1-x)}\right]_2=0\ .
			\end{equation}
			Une forme plus agréable de cette relation est
			\begin{equation}\label{eq: five term birapport}
			\sum_{i=0}^4 (-1)^i [r(a_0,\ldots,\widehat{a_i},\ldots,a_4)]_2=0\ ,
			\end{equation}
			où $a_0,a_1,a_2,a_3,a_4$ sont des points deux à deux distincts de la droite projective $\mathbb{P}^1(F)$ et où l'on note $r(a,b,c,d)=\frac{(a-c)(b-d)}{(a-d)(b-c)}$ le birapport de $4$ points. On obtient~\eqref{eq: five term abstract} en spécialisant à $(a_0,a_1,a_2,a_3,a_4)=(\infty, 0,1,x,y)$, et les relations $[1]_2=0$, $[x]_2+[x^{-1}]_2=0$, et~$[x]_2+[1-x]_2=0$ en spécialisant à des cas dégénérés. Les travaux de \textcite{suslinbloch} impliquent\footnote{Voir la remarque \ref{rema: suslin bloch} plus bas. Dans le cas d'un corps $F$ algébriquement clos, un argument plus explicite dû à \textcite{wojtkowiakfunctional} donne le même résultat. Cet argument montre aussi que les relations à~$5$ termes engendrent les relations fonctionnelles (linéaires et en une variable) de la fonction $\mathrm{P}_2$. On renvoie aux travaux récents de \textcite{dejeufunctional} pour le cas des relations en plusieurs variables.} que ces relations engendrent les relations polylogarithmiques pour~$n=2$, quel que soit le corps $F$. 
			
			La preuve de la conjecture de Zagier pour $n=3$ repose notamment sur une variante de la relation~\eqref{eq: five term birapport} découverte par \textcite{goncharovseattle, goncharovgeometry}. On peut l'écrire sous la forme très (anti)symétrique, découverte indépendamment par Goncharov et par Zagier, d'une relation \og à 840 termes \fg{} dans $\mathcal{B}_3(F)$ :
			\begin{equation}\label{eq: trirapport}
			\sum_{i=0}^6 (-1)^i R_3(a_0,\ldots,\widehat{a_i},\ldots,a_6) = 0\ ,
			\end{equation}
			où $a_0,\ldots,a_6$ sont des points de $\mathbb{P}^2(F)$ en position générale et $R_3(a_1,\ldots,a_6)\in \QQ[F^\times]$ est un invariant d'une configuration de $6$ points en position générale dans le plan projectif obtenue en antisymétrisant une variante du birapport (appelée \emph{trirapport}). Goncharov conjecture que les relations~\eqref{eq: trirapport} et leurs spécialisations engendrent toutes les relations polylogarithmiques pour $n=3$.

	\subsection{Complexes polylogarithmiques}
	
		Le programme de Goncharov va bien au-delà du cas des corps de nombres et propose une description \og symbolique \fg{} de la $K$-théorie rationnelle d'un corps $F$ quelconque, dans l'esprit de la description de $K_2(F)$ par \textcite{matsumoto} et de sa généralisation par \textcite{milnorktheory}. Cette proposition, qui à première vue peut sembler arbitraire, sera éclairée par des considérations motiviques dans le prochain paragraphe.

		Les morphismes $\widetilde{\delta}_n$ définis par \eqref{eq: delta deux tilde} et \eqref{eq: delta tilde} passent au quotient et induisent des morphismes 
		$$\delta_2:\mathcal{B}_2(F)\to \Lambda^2F^\times_\QQ \;\;\hbox{ et }\;\;\delta_n:\mathcal{B}_n(F)\to \mathcal{B}_{n-1}(F)\otimes F^\times_\QQ \;\hbox{ pour }\; n\geq 3\ .$$ 
		Ces morphismes s'assemblent en un complexe, le $n$-ième \emph{complexe polylogarithmique}, noté $\mathcal{B}^\bullet(F,n)$, où $\mathcal{B}_n(F)$ est placé en degré cohomologique $1$ :
		\begin{equation*}
		\begin{split}
		0\to \mathcal{B}_n(F)\To \mathcal{B}_{n-1}(F)\otimes F^\times_\QQ \To \mathcal{B}_{n-2}(F&)\otimes \Lambda^2F^\times_\QQ \To \cdots \\
		& \cdots \To \mathcal{B}_2(F)\otimes \Lambda^{n-2}F^\times_\QQ \To \Lambda^nF^\times_\QQ\to 0\ .
		\end{split}
		\end{equation*}
		Il est commode d'étendre les définitions à $n\leq 1$ en posant $\mathcal{B}_0(F)=\QQ$, $\mathcal{B}_1(F)=F^\times_\QQ$, et donc $\mathcal{B}^\bullet(F,0)=\QQ$, $\mathcal{B}^\bullet(F,1)=F^\times_\QQ[-1]$. 
		
		On peut maintenant énoncer une conjecture centrale du programme de Goncharov. Notons~$\mathrm{gr}_\gamma$ les gradués pour la $\gamma$-filtration en $K$-théorie \parencite{souleoperations}, issue de la structure de $\lambda$-anneau sur~$K_\bullet(F)_\QQ$ induite par les opérateurs de puissance extérieure des $F$-espaces vectoriels.
		
		\begin{conj}\label{conj: goncharov polylog complexes} On a des isomorphismes fonctoriels en le corps $F$, pour tout entier $n\geq 0$ et tout entier relatif $i$ :
		$$H^i(\mathcal{B}^\bullet(F,n)) \simeq \mathrm{gr}^n_\gamma K_{2n-i}(F)_\QQ\ .$$
		\end{conj}
		
		On voit facilement que cette conjecture est vérifiée pour $n\leq 1$, ainsi que dans le cas $i=n$ puisque $H^n(\mathcal{B}^\bullet(F,n))= K^{\mathrm{M}}_n(F)_\QQ$, la $K$-théorie de Milnor rationnelle, qui est isomorphe à $\mathrm{gr}^n_\gamma K_n(F)_\QQ$ d'après \textcite{suslinhomology} et \textcite{souleoperations}. 
		
		\begin{rema}\label{rema: suslin bloch}
		Une des sources de la conjecture \ref{conj: goncharov polylog complexes} est le résultat de \textcite{suslinbloch} qui fournit une suite exacte\footnote{Le résultat de Suslin est plus précis puisqu'il prend aussi en compte la torsion en $K$-théorie. Suslin ne traite que le cas $F$ infini, le cas $|F|\geq 4$ ayant été traité par \textcite[VI, Theorem 5.2]{weibelkbook}. Les cas $|F|<4$ ne posent aucun problème si l'on néglige la torsion.}
		$$0\To \mathrm{gr}^2_\gamma K_3(F)_\QQ \To \QQ[F^\times]/\mathcal{R}_2^{(5)}(F) \stackrel{\delta_2}{\To} \Lambda^2(F^\times_\QQ) \To K_2(F)_\QQ\To 0\ ,$$
		où l'on a noté $\mathcal{R}_2^{(5)}(F)\subset \mathcal{R}_2(F)$ le sous-espace engendré par les relations à $5$ termes~\eqref{eq: five term abstract} et leurs spécialisations. Combiné à l'isomorphisme $\mathrm{gr}^2_\gamma K_3(F)\stackrel{\sim}{\to}\mathrm{gr}^2_\gamma K_3(F(t))$, cette suite exacte montre qu'on a l'égalité $\mathcal{R}_2^{(5)}(F)=\mathcal{R}_2(F)$. On en déduit un isomorphisme~$H^1(\mathcal{B}^\bullet(F,2))\simeq \mathrm{gr}^2_\gamma K_3(F)_\QQ$.
		\end{rema}
		
		Dans le cas d'un corps de nombres, les travaux de Borel (et la réinterprétation du régulateur de Borel par Beilinson, voir la remarque \ref{rema: regulateurs}) ont pour conséquence l'annulation :
		\begin{equation}\label{eq: annulation K theorie corps de nombres}
		\mathrm{gr}^n_\gamma K_{2n-i}(F)_\QQ = 0  \;\;\;\mbox{ si } \; i\neq 1 \;\;\;\;\; (F \mbox{ corps de nombres}) \ .
		\end{equation}
		 La conjecture~\ref{conj: goncharov polylog complexes} implique donc dans ce cas un isomorphisme $K_{2n-1}(F)_\QQ\simeq \ker(\delta_n)$ dont la composition avec l'inclusion $\ker(\delta_n)\subset \mathcal{B}_n(F)$ devrait donner lieu à un morphisme $\varphi_n$ comme dans la conjecture~\ref{conj: zagier k theorie}.
		
		Notons aussi que la conjecture~\ref{conj: goncharov polylog complexes} est compatible avec (et en fait implique) la conjecture d'annulation suivante de \textcite{beilinsonhigher, souleoperations}.
		
		\begin{conj}[Beilinson--Soulé]\label{conj: beilinsonsoule} Pour tout corps $F$ et tout entier $n\geq 1$ :
		$$\mathrm{gr}^n_\gamma K_{2n-i}(F)_\QQ=0 \;\;\; \hbox{ si }\; i\leq 0\ .$$
		\end{conj}
		
		Cette conjecture est vérifiée pour les corps de nombres par~\eqref{eq: annulation K theorie corps de nombres}.
		
		\begin{rema}\label{rema: complexe motivique}
		La conjecture~\ref{conj: goncharov polylog complexes} exprime le fait que $\mathcal{B}^\bullet(F,n)$ est un \emph{complexe motivique} pour la variété $\mathrm{Spec}(F)$, au sens de \textcite{lichtenbaumnonnegative} et de \textcite{beilinsonheight}.
		\end{rema}

		\subsection{Une stratégie vers la conjecture de Zagier, première version}\label{par: strategie premiere}
					
			En extrapolant à partir des preuves de la conjecture de Zagier pour $n=2,3$, il est tentant, en suivant \textcite{goncharovexplicit, goncharovseattle, goncharovgeometry}, de procéder comme suit. On cherche d'abord une application~$R_n$ de l'ensemble des familles de~$2n$ points dans~$\mathbb{P}^{n-1}(F)$ en position générale, à valeurs dans $\QQ[F^\times]$. Supposons que $R_n$~soit invariante pour l'action de~$\mathrm{GL}_n(F)$ sur~$\mathbb{P}^{n-1}(F)$ et qu'on ait la relation dans $\mathcal{B}_n(F)$ :
			\begin{equation}\label{eq: equation fonctionnelle generale}
			\sum_{i=0}^{2n}(-1)^i R_n(a_0,\ldots,\widehat{a_i},\ldots,a_{2n})=0\ ,
			\end{equation}
			pour tous les $a_0,\ldots,a_{2n}$ en position générale dans $\mathbb{P}^{n-1}(F)$. (On peut même aller plus loin et spéculer sur le fait que les relations~\eqref{eq: equation fonctionnelle generale} et leurs spécialisations engendrent les relations polylogarithmiques.)
			
			La suite de l'argument repose sur une construction d'algèbre homologique due à  \textcite{suslinhomology}. Notons $G_k(F,r)$ le $\QQ$-espace vectoriel librement engendré par l'ensemble des $\mathrm{GL}_r(F)$-orbites de familles $(a_0,\ldots,a_k)$ de points de $F^r$ en position générale. On a des différentielles
			$$\partial : G_k(F,r)\To G_{k-1}(F,r) \; ,\; (a_0,\ldots,a_k)\mapsto \sum_{i=0}^k(-1)^i(a_0,\ldots,\widehat{a_i},\ldots,a_k)$$
			$$\partial': G_k(F,r) \To G_{k-1}(F,r-1) \; , \; (a_0,\ldots,a_k)\mapsto \sum_{i=0}^k(-1)^i(a_i | a_0,\ldots,\widehat{a_i},\ldots,a_k)\ ,$$
			où la notation $(a_i | a_0,\ldots,\widehat{a_i},\ldots,a_k)$ désigne la configuration associée aux images des points $a_0,\ldots,\widehat{a_i},\ldots,a_k$ dans $F^r/Fa_i\simeq F^{r-1}$. Elles donnent lieu à un complexe double formé des $G_k(F,r)$ pour $k\geq r\geq n$, dont le complexe total est noté $BG^\bullet(F,n)$, le~$n$\nobreakdash-ième \emph{complexe bigrassmannien}. Par convention, $G_k(F,r)$ y est placé en degré cohomologique~$2n-k$. La motivation pour cette construction est que, si $F$ est un corps infini, le complexe~$(G_\bullet(F,n),\partial)$ est obtenu en prenant les coinvariants d'une résolution du~$\mathrm{GL}_n(F)$-module trivial $\QQ$, et permet donc de construire des classes de cohomologie de $\mathrm{GL}_n(F)$. On considère tous les $n$ à la fois pour produire des classes de cohomologie du groupe linéaire stable\footnote{La conjecture du rang pour les corps de nombres, qui est connue grâce aux travaux de \textcite{yang} complétés par \textcite{borelyang}, permet de s'affranchir de cette étape.} $\mathrm{GL}(F)$.
			
			Une application $R_n$ comme ci-dessus induit un morphisme $R_n:G_{2n-1}(F,n)\to \mathcal{B}_{n}(F)$ qui vérifie $R_n\circ \partial=0$. Supposons qu'on puisse l'étendre en un morphisme de complexes 
			$$\widetilde{R}_n: BG^\bullet(F,n)\To \mathcal{B}^\bullet(F,n)\ .$$
			Alors on obtient par un argument d'homologie des groupes un morphisme
			$$H_{2n-1}(\mathrm{GL}(F),\QQ) \To \ker(\delta_n)\subset \mathcal{B}_n(F)\ ,$$
			et donc un morphisme $\varphi_n:K_{2n-1}(F)_\QQ\to \mathcal{B}_n(F)$. Il n'y a aucune raison \emph{a priori} pour que ce morphisme fasse commuter le diagramme \eqref{eq: diag comm zagier k theorie} (ou même pour qu'il soit non nul). Des techniques générales développées par Goncharov permettent de le montrer si l'on sait calculer la composante $G_n(F,n)\to \Lambda^n(F^\times)_\QQ$ de $\widetilde{R}_n$, à l'autre extrémité du complexe bigrassmannien.
			
		\subsection{Nécessité des polylogarithmes multiples}			
			
			L'existence d'un morphisme $\widetilde{R}_n$ comme dans le paragraphe précédent est appelée \og conjecture optimiste \fg{} par \textcite[Conjecture 5.12]{goncharovseattle}. En effet, cette stratégie semble être difficile à faire fonctionner au-delà des cas $n=2,3$. La raison est qu'il est difficile de produire des relations \og intéressantes \fg{} dans $\mathcal{B}_n(F)$ (ou, ce qui revient au même, des équations fonctionnelles intéressantes pour les polylogarithmes classiques $\mathrm{Li}_n$ ou $\mathrm{P}_n$) pour $n\geq 4$, et les relations de cocycle~\eqref{eq: equation fonctionnelle generale} semblent inaccessibles en général. En pratique, les équations fonctionnelles pour les polylogarithmes classiques font naturellement intervenir les polylogarithmes \emph{multiples}, dont on ne peut se passer pour $n\geq 4$. Il faut donc adapter la stratégie en considérant cette plus grande classe de fonctions  --- même si au bout du compte seuls les polylogarithmes \emph{classiques} sont censés intervenir si l'on croit à la conjecture de Zagier. 
			
			Dit autrement, il apparaît rétrospectivement que les espaces~$\mathcal{B}_n(F)$ et les complexes~$\mathcal{B}^\bullet(F,n)$ sont des objets trop petits pour être facilement maniables. On les verra plus bas à l'intérieur d'objets plus naturels et plus structurés : la cogèbre de Lie motivique~$\mathcal{C}(F)$ et son complexe de Chevalley--Eilenberg. Le prochain paragraphe a pour but d'introduire ces objets et les liens qu'ils entretiennent avec la $K$-théorie.

\section{Motifs de Tate mixtes et $K$-théorie}

On introduit maintenant la catégorie des motifs de Tate mixtes sur $F$ et on développe son formalisme tannakien. Il sera commode de supposer que la conjecture d'annulation de Beilinson--Soulé (conjecture~\ref{conj: beilinsonsoule}) est vraie, ou, ce qui revient au même, de se restreindre au cas des corps $F$ pour laquelle elle est vérifiée, ce qui inclut le cas des corps de nombres.

	\subsection{Motifs de Tate mixtes}
	
		Pour un corps $F$, notons $\DM(F)$ la catégorie triangulée des motifs sur $F$ à coefficients dans $\QQ$ (dont différentes constructions équivalentes existent, notamment par \textcite{hanamuramixed}, \textcite{levinemixed}, et \textcite{voevodskytriangulated}). Chaque variété $X$ sur $F$ définit un objet de cette catégorie, le \emph{motif} de $X$, qu'il faut voir comme un complexe qui contrôle toutes les versions de la cohomologie de $X$ (de Rham, singulière, $\ell$-adique, etc.). Des foncteurs de réalisation définis sur $\DM(F)$ permettent de retrouver ces différents groupes de cohomologie à partir du motif de $X$. La différence fondamentale entre la catégorie des motifs et les différentes catégories de réalisations est que les morphismes dans $\DM(F)$ sont définis en termes de cycles algébriques et sont notamment reliés aux groupes de Chow supérieurs de \textcite{blochhigher}.
		
		 La catégorie $\DM(F)$ est tensorielle et on dispose d'un objet $\otimes$-inversible $\QQ(-1)$, le \emph{motif de Lefschetz}, dont les réalisations sont les groupes de cohomologie $H^1(\mathbb{P}^1_F\setminus \{0,\infty\})$. On note $\QQ(-n)=\QQ(-1)^{\otimes n}$ pour tout $n\in\ZZ$. Grâce à la comparaison entre groupes de Chow supérieurs et $K$-théorie on a des isomorphismes 
		\begin{equation}\label{eq: hom DM K theorie}
		\mathrm{Hom}_{\DM(F)}(\QQ(-n),\QQ(0)[i]) \simeq \mathrm{gr}^n_\gamma K_{2n-i}(F)_\QQ\ ,
		\end{equation}
		qui expliquent l'irruption de la théorie des motifs dans des questions liées à la $K$-théorie.
		
		Grâce à l'annulation qui découle de la conjecture de Beilinson--Soulé, \textcite{levinetate} définit une $t$-structure sur la sous-catégorie triangulée de~$\DM(F)$ engendrée par les~$\QQ(-n)$ pour~$n\in\ZZ$, dont le cœur est par définition la catégorie~$\MT(F)$ des \emph{motifs de Tate mixtes} sur~$F$. Tout objet~$M$ de cette catégorie est muni d'une filtration finie (par le \emph{poids}) par des sous-objets~$W_{2n}M$ de telle sorte que chaque gradué~$\mathrm{gr}_{2n}^WM$ est une somme directe finie d'objets~$\QQ(-n)$. 
		
		En suivant \textcite{beilinsonheight} on est amené à formuler la conjecture additionnelle suivante, parfois appelée \emph{conjecture du $K(\pi,1)$}.
		
		\begin{conj}\label{conj: beilinson ext}
		Pour tous entiers $n,i$ le morphisme naturel
		$$\mathrm{Ext}^i_{\MT(F)}(\QQ(-n),\QQ(0)) \To \mathrm{Hom}_{\DM(F)}(\QQ(-n),\QQ(0)[i])\simeq \mathrm{gr}^n_\gamma K_{2n-i}(F)_\QQ$$
		est un isomorphisme.
		\end{conj}
		
		Le morphisme ci-dessus est, pour des raisons générales, un isomorphisme si~$i=0,1$ et injectif pour $i=2$. On en déduit, grâce à l'annulation \eqref{eq: annulation K theorie corps de nombres}, que la conjecture est vérifiée dans le cas d'un corps de nombres.
		
		Il est facile de décrire géométriquement les extensions de~$\QQ(-1)$ par~$\QQ(0)$ dans~$\MT(F)$. Via l'isomorphisme entre le groupe d'extensions correspondant et~$\mathrm{gr}_\gamma^1K_1(F)_\QQ=F^\times_\QQ$, un élément $x\in F^\times$ est représenté par le \emph{motif de Kummer}~$K(x)$, dont les réalisations sont les groupes de cohomologie relative $H^1(\mathbb{P}^1_F\setminus \{0,\infty\},\{1,x\})$ si $x\neq 1$. La question de la description géométrique des extensions de~$\QQ(-n)$ par~$\QQ(0)$ dans~$\MT(F)$ pour~$n\geq 2$ est difficile, et constitue en quelque sorte le coeur de la conjecture de Zagier, comme nous allons le voir.
		
	\subsection{Formalisme tannakien et complexes motiviques}
		
		La catégorie $\MT(F)$ est tannakienne neutre à coefficients dans $\QQ$. Les différentes réalisations donnent lieu à des foncteurs fibre, parmi lesquels les foncteurs fibre de Rham (si $F$ est de caractéristique nulle) $\omega_{\mathrm{dR}}:\MT(F)\to \mathrm{Vect}_F$, et Betti (cohomologie singulière, relative à un plongement complexe $\sigma:F\to \CC$) $\omega_{\mathrm{B}}:\MT(F)\to \mathrm{Vect}_\QQ$. La comparaison $\omega_{\mathrm{dR}}\otimes_F\CC \stackrel{\sim}{\to} \omega_{\mathrm{B}}\otimes_\QQ\CC$ entre ces foncteurs donne lieu à des \emph{matrices des périodes} pour les objets de $\MT(F)$. La matrice des périodes de $\QQ(-n)$ est $\left((2\mathrm{i}\pi)^n\right)$. 
		
		Dans la suite on travaillera surtout avec le foncteur fibre $\mathbb{Z}$-gradué canonique à valeurs dans les $\QQ$-espaces vectoriels :
		\begin{equation}\label{eq: foncteur fibre canonique}
		\omega: M \mapsto\bigoplus_{n\in\mathbb{Z}} \mathrm{Hom}_{\MT(F)}(\QQ(-n),\mathrm{gr}_{2n}^WM)\ .
		\end{equation}
		Si $F$ est de caractéristique nulle, il s'agit d'une structure rationnelle sur le foncteur fibre de Rham : $\omega\otimes_\QQ F\simeq \omega_{\mathrm{dR}}$. La graduation sur $\omega$ est couramment appelée \emph{poids}, de sorte que $\QQ(-n)$ est de poids $n$\footnote{Les poids sont donc divisés par $2$ par rapport à la notion habituelle, pour laquelle $\QQ(-n)$ est de poids $2n$.}.
		
		Par le formalisme tannakien, la catégorie $\MT(F)$ est équivalente à la catégorie des comodules sur une algèbre de Hopf commutative $\mathbb{Z}$-graduée $\mathcal{H}(F)$, où la graduation est encore appelée \emph{poids}. Concrètement, $\mathcal{H}_n(F)$ est engendré en tant que $\QQ$-espace vectoriel par des coefficients matriciels $(M,v,\varphi)$, où $M$ est un objet de $\MT(F)$ et où~$v:\QQ(-n)\to \mathrm{gr}_{2n}^WM$ et $\varphi:\mathrm{gr}_0^WM\to \QQ(0)$ sont deux morphismes. Les relations entre ces coefficients matriciels sont engendrées par les égalités $(M,f\circ v,\psi)=(N,v,\psi\circ f)$ pour des morphismes $f:M\to N$ et $v:\QQ(-n)\to \mathrm{gr}_{2n}^WM$, $\psi:\mathrm{gr}_0^WN\to \QQ(0)$. On a  $\mathcal{H}_n(F)=0$ si $n<0$, $\mathcal{H}_0(F)=\QQ$, et $\mathcal{H}_1(F)=F^\times_\QQ$.  
		
		Le produit dans $\mathcal{H}(F)$ est induit par le produit tensoriel des objets de $\MT(F)$ et le coproduit~$\Delta:\mathcal{H}(F)\to \mathcal{H}(F)\otimes \mathcal{H}(F)$ se calcule par la formule
		$$\Delta(M,v,\varphi) = \sum_i (M,e_i,\varphi)\otimes (M,v,e_i^\vee)\ ,$$
		où $(e_i)$ est n'importe quelle base de $\omega(M)$.
		
		Il sera plus pratique de travailler modulo produits et de manipuler la cogèbre de Lie des indécomposables de $\mathcal{H}(F)$
		$$\mathcal{C}(F) = \mathcal{H}_{>0}(F)/\mathcal{H}_{>0}(F)\mathcal{H}_{>0}(F) \;\; ,\;\; \delta: \mathcal{C}(F)\to \Lambda^2\mathcal{C}(F)\ ,$$
		dont le \emph{cocrochet} $\delta$ est induit par $\Delta$. Elle est $\mathbb{Z}_{\geq 1}$-graduée et vérifie $\mathcal{C}_1(F)=F^\times_\QQ$. Pour $x\in F^\times$, l'élément correspondant de $\mathcal{C}_1(F)$ est le coefficient matriciel du motif de Kummer $K(x)$ où $v$ et $\varphi$ sont induits par les isomorphismes naturels $\mathrm{gr}_2^WK(x)\simeq \QQ(-1)$ et $\mathrm{gr}_0^W K(x)\simeq \QQ(0)$ respectivement. On le note
		$$\log^{\mathcal{C}}(x) \in \mathcal{C}_1(F)\ .$$
		Les éléments de $\mathcal{H}(F)$ ou $\mathcal{C}(F)$ sont parfois appelés \emph{périodes motiviques}, un terme qui recouvre plusieurs variantes de l'idée des coefficients matriciels dans une catégorie tannakienne de motifs (voir la remarque~\ref{rema: periodes motiviques}).
		
		\begin{rema}
		Le groupe tannakien correspondant au foncteur fibre~\eqref{eq:
                  foncteur fibre canonique}, aussi appelé \emph{groupe de Galois motivique} de la catégorie $\MT(F)$, est le produit semi-direct de~$\mathbb{G}_m$ avec le schéma en groupes pro-unipotent
                $U(F)=\mathrm{Spec}(\mathcal{H}(F))$. Il
                est complètement déterminé par l'algèbre de Lie graduée
                pro-nilpotente~$\mathfrak{u}(F)$ de~$U(F)$ (via la formule de
                Baker--Campbell--Hausdorff), qui s'obtient en dualisant la cogèbre de Lie~$\mathcal{C}(F)$. L'avantage de travailler dans $\mathcal{H}(F)$ ou $\mathcal{C}(F)$ est qu'on peut facilement en produire des éléments et calculer leur coproduit ou cocrochet, alors qu'exhiber des points de~$U(F)$ ou des éléments de~$\mathfrak{u}(F)$ est incommode.
		\end{rema}

		Pour un entier $n\geq 1$ on considère la partie de poids $n$ du complexe de Chevalley--Eilenberg de $\mathcal{C}(F)$, noté $\mathrm{CE}^\bullet(\mathcal{C}(F))_n$, où $\mathcal{C}_n(F)$ est placé en degré cohomologique $1$ :
		$$0\to \mathcal{C}_n(F) \To (\Lambda^2\mathcal{C}(F))_n \To (\Lambda^3\mathcal{C}(F))_n \To \cdots \To (\Lambda^{n-1}\mathcal{C}(F))_n \To  \Lambda^n\mathcal{C}_1(F)\to 0\ .$$
		Sa cohomologie calcule des groupes d'extensions dans la catégorie des $\mathcal{C}(F)$-comodules gradués, et plus précisément :
		\begin{equation}\label{eq: isom H CE Ext}
		H^i(\mathrm{CE}^\bullet(\mathcal{C}(F))_n) \simeq \mathrm{Ext}^i_{\MT(F)}(\QQ(-n),\QQ(0))\ .
		\end{equation}

		Si l'on croit à la conjecture~\ref{conj: beilinson ext}, les
                complexes $\mathrm{CE}^\bullet(\mathcal{C}(F))_n$ sont donc
                des \emph{complexes motiviques} (voir la remarque~\ref{rema: complexe motivique}), dont l'existence est
                conditionnelle à la conjecture de Beilinson--Soulé. (On a vu que les deux
                conjectures sont vérifiées dans le cas d'un corps de nombres.)
                Une partie des travaux de Goncharov vise à définir, en
                s'inspirant de $\mathrm{CE}^\bullet(\mathcal{C}(F))$, des
                complexes motiviques inconditionnels qui soient aussi
                explicites et \og petits \fg{} que possible.
		
		Notons que dans le cas d'un corps de nombres $F$ on obtient l'isomorphisme inconditionnel suivant qui identifie la $K$-théorie rationnelle de $F$ avec l'espace des primitifs de la cogèbre de Lie motivique :
		\begin{equation}\label{eq: isom K theorie primitifs C}
		K_{2n-1}(F)_\QQ \simeq H^1(\mathrm{CE}^\bullet(\mathcal{C}(F))_n) = \ker\big(\delta:\mathcal{C}_n(F)\to (\Lambda^2\mathcal{C}(F))_n\big)\ .
		\end{equation}
		De plus, on voit grâce à \eqref{eq: annulation K theorie corps de nombres}, \eqref{eq: isom H CE Ext}, et le fait que la conjecture \ref{conj: beilinson ext} est vérifiée pour les corps de nombres, que $\mathcal{C}(F)$ est colibre. Dualement, l'algèbre de Lie motivique $\mathfrak{u}(F)$ est libre et ses indécomposables en poids $-n\leq -1$ sont donnés par l'espace $K_{2n-1}(F)^\vee_\QQ$.
		
		Dans le paragraphe suivant on explique comment penser au régulateur dans le cadre du formalisme tannakien.
		
	\subsection{Réalisation de Hodge et régulateur}
		
		Soit $F$ un corps de nombres. On a pour tout plongement $\sigma:F\to \CC$ un foncteur de réalisation de Hodge $\MT(F)\to \MHS$ vers la catégorie des structures de Hodge mixtes rationnelles. 
		Il induit pour tout~$n\geq 1$ un morphisme au niveau des groupes d'extensions
		\begin{equation}\label{eq: hodge real ext}
		\mathrm{Ext}^1_{\MT(F)}(\QQ(-n),\QQ(0)) \To \mathrm{Ext}^1_{\MHS}(\QQ(-n),\QQ(0)) \simeq \CC/(2\mathrm{i}\pi)^n\QQ\ .
		\end{equation}
		Rappelons qu'un nombre complexe $a\in\CC/(2\mathrm{i}\pi)^n\QQ$ correspond par ce dernier isomorphisme à la structure de Hodge mixte dont une matrice des périodes est 
		$$\left(\begin{matrix}1 & a \\ 0 & (2\mathrm{i}\pi)^n\end{matrix}\right)\ ,$$
		où la deuxième colonne donne la position de la filtration de Hodge dans une base rationnelle compatible à la filtration par le poids. Le motif de Kummer $K(x)$ donne lieu à une structure de Hodge mixte pour laquelle 
		\begin{equation}\label{eq: calcul log}
		a=\int_1^{z}\frac{dt}{t} = \log(z) \;\;\in \CC/2\mathrm{i}\pi\QQ \ ,
		\end{equation}
		avec $z=\sigma(x)$.
		
		Si l'on somme~\eqref{eq: hodge real ext} sur tous les plongements complexes de $F$, on obtient en composant avec $\CC/(2\mathrm{i}\pi)^n\QQ\to \CC/(2\mathrm{i}\pi)^n\RR$ et en utilisant~\eqref{eq: identification Rdn} et~\eqref{eq: hom DM K theorie} un morphisme
		\begin{equation}\label{eq: regulateur hodge}
		\varpi_n:K_{2n-1}(F)_\QQ\To \RR^{d_n}
		\end{equation}
		qu'on peut appeler \emph{régulateur de Hodge}. On devrait avoir l'égalité $\varpi_n=\rho_n$ avec le régulateur de Borel, mais ce résultat ne semble pas être dans la littérature.

		Le point de vue de la réalisation de Hodge fournit en fait plus que le régulateur~\eqref{eq: regulateur hodge} : on a pour tout entier $n$ un morphisme
		$$\varpi_n:\mathcal{C}_n(F)\To \RR^{d_n}$$
		qui induit~\eqref{eq: regulateur hodge} en restriction au noyau du cocrochet $\delta$, isomorphe à $K_{2n-1}(F)_\QQ$ via~\eqref{eq: isom K theorie primitifs C}. Il apparaît chez \textcite{goncharovvolumes} sous le nom de morphisme des \emph{périodes réelles} ou \emph{Lie-périodes} (voir aussi \textcite{brownsvmzv} pour une interprétation tannakienne comme variante du morphisme des \emph{périodes univaluées}). Explicitons sa définition en suivant \textcite{beilinsondeligne}. Représentons un élément de $\mathcal{C}_n(F)$ par un coefficient matriciel $(M,v,\varphi)$. Pour un plongement $\sigma:F\to \CC$, soit $P$ une matrice des périodes de~$M$ (matrice de l'isomorphisme naturel $\omega(M)\otimes_\QQ\CC\simeq \omega_{\mathrm{dR}}(M)\otimes_{F}\CC \stackrel{\sim}{\to} \omega_{\mathrm{B}}(M)\otimes_\QQ\CC$). 
		Si~$D$ désigne la matrice diagonale qui agit par $(-1)^k$ en poids $k$, le produit $D\overline{P}^{-1}P$ est alors une matrice unipotente. La composante de $\varpi_n(M,v,\varphi)$ correspondant au plongement~$\sigma$ est alors l'accouplement
		$$\langle \, \varphi \, , \, \textstyle\frac{1}{2}\log(D\overline{P}^{-1}P)v  \, \rangle \ .$$
		Par exemple, $\varpi_1(\log^{\mathcal{C}}(x))$ a pour composantes les $\log|\sigma(x)|$.
		
\section{La conjecture de Zagier dans le contexte motivique}

	On passe maintenant à l'interprétation motivique de la conjecture de Zagier, d'après \textcite{beilinsondeligne} et \textcite{goncharovseattle, goncharovgeometry}. En suivant Goncharov, on arrivera naturellement à l'énoncé de la \emph{conjecture de liberté} sur la structure de la cogèbre de Lie motivique.
	
	\subsection{Motifs polylogarithmiques et polylogarithmes motiviques}
	
		On dispose pour tout $x\in F^\times$ d'un \emph{motif polylogarithmique} $L(x)$, qui est un ind-objet de $\MT(F)$
		 tel que $\mathrm{gr}_{2k}^W L(x)\simeq \QQ(-k)$ pour tout entier naturel $k$. Une matrice des périodes de $W_{2n}L(x)$, relative à un plongement $\sigma:F\to \CC$, est de la forme suivante, où l'on note $z=\sigma(x)$.
		\begin{equation}\label{eq: matrice periodes Li}
		\left(\begin{matrix}
		1 & \Li_1(z) & \Li_2(z)        & \Li_3(z) & \cdots & \Li_n(z)\\
		 & 2\mathrm{i}\pi    & 2\mathrm{i}\pi \log(z)  & 2\mathrm{i}\pi \frac{\log^2(z)}{2}& \cdots & 2\mathrm{i}\pi \frac{\log^{n-1}(z)}{(n-1)!}\\
		 &            & (2\mathrm{i}\pi)^2      & (2\mathrm{i}\pi)^2 \log(z) & \cdots & (2\mathrm{i}\pi)^2 \frac{\log^{n-2}(z)}{(n-2)!} \\
		 &  &  & (2\mathrm{i}\pi)^3 & &  \\
		 & & & & & \vdots \\
		& & 0 & & \ddots &  \\
		& & & & & \\
		& & & & &  (2\mathrm{i}\pi)^n 
		\end{matrix}\right)
		\end{equation}
		(Si $x=1$, $\Li_1(z)$ n'est pas défini et doit être remplacé par $0$ dans cette matrice.) L'apparition de $\Li_n(z)$ dans une matrice des périodes s'explique par la formule intégrale suivante, valable pour $z\notin \, ]1,+\infty[$ :
		$$\Li_n(z) = \int_{0\leq t_1\leq\cdots \leq t_n\leq 1} \frac{z\, dt_1\cdots dt_n}{(1-zt_1)t_2\cdots t_n} \ \cdot$$
		
		Pour tout entier $n\geq 1$, le motif polylogarithmique donne lieu, via les isomorphismes $\mathrm{gr}^W_{0}L(x)\simeq \QQ(0)$ et $\mathrm{gr}^W_{2n}L(x)\simeq \QQ(-n)$, à un coefficient matriciel noté
		\begin{equation}\label{eq: polylog motivique}
		\Li_n^{\mathcal{C}}(x) \in \mathcal{C}_n(F)
		\end{equation}
		et appelé \emph{polylogarithme motivique}. 
		La structure très simple de $L(x)$, reflétée dans la structure par blocs de la matrice~\eqref{eq: matrice periodes Li}, permet de calculer son cocrochet, donné pour $n\geq 2$ par la formule :
		\begin{equation}\label{eq: cocrochet Li n}
		\delta\Li_n^{\mathcal{C}}(x) = \Li_{n-1}^{\mathcal{C}}(x) \wedge \log^{\mathcal{C}}(x)\ .
		\end{equation}
		Notons que $\delta\Li_n^{\mathcal{C}}(x)$ a seulement une composante non nulle dans $\mathcal{C}_{n-1}(F)\wedge \mathcal{C}_1(F)$. Si~$F$ est un corps de nombres, un calcul rapide utilisant~\eqref{eq: matrice periodes Li} montre qu'on a l'égalité au niveau des régulateurs de Hodge :
		\begin{equation}\label{eq: periode reelle Li n}
		\varpi_n(\mathrm{Li}_n^{\mathcal{C}}(x)) = \mathrm{P}_n^F(x)\ .
		\end{equation} 
		Cela justifie le choix de $\mathrm{P}_n$ comme version univaluée de $\Li_n$.
	
\begin{rema}\label{rema: periodes motiviques}
  C'est l'équation~\eqref{eq: periode reelle Li n} qui justifie l'appellation
  \emph{polylogarithme motivique} pour l'élément~$\Li_n^{\mathcal{C}}(x)$. Il
  existe dans la littérature des variantes de cette construction qui portent
  le même nom et jouent des rôles différents, mais qui ont en commun de venir
  de coefficients matriciels dans une catégorie tannakienne de motifs (comme~$\MT(F)$) ou de systèmes de réalisations (structures de Hodge mixtes,
  etc.). Par exemple, on peut considérer une version~$\Li_n^{\mathcal{H}}(x)$
  dans l'algèbre de Hopf motivique~$\mathcal{H}(F)$, dont
  $\Li_n^{\mathcal{C}}(x)$~est la version modulo produits. On en déduit par
  itération maximale du coproduit (réduit) le \emph{symbole} de~$\Li_n(x)$,
  qui est égal à
  $\Li_n^{\mathcal{S}}(x)=-(1-x)\otimes x\otimes\cdots \otimes x$ et vit dans
  $(\mathcal{H}_1(F))^{\otimes n} = (F^\times_\QQ)^{\otimes n}$. Une autre
  variante, plus proche de la notion de période, est de considérer l'algèbre
  des périodes motiviques, relative à un plongement complexe de~$F$, notée
$\mathcal{P}(F)=\mathcal{O}(\underline{\mathrm{Isom}}^{\otimes}_{\MT(F)}(\omega_{\mathrm{dR}},\omega_{\mathrm{B}}))$, dont les éléments sont représentés par des coefficients matriciels $(M,\alpha,\tau)$ avec~$M\in \MT(F)$, $\alpha\in \omega_{\mathrm{dR}}(M)$ et $\tau\in \omega_{\mathrm{B}}(M)^\vee$. Cette algèbre est munie d'un morphisme d'algèbres vers~$\CC$, le \emph{morphisme des périodes}, et on peut relever (une détermination de) $\Li_n(x)\in \CC$ en un élément $\Li_n^{\mathcal{P}}(x)\in \mathcal{P}(F)$. Cette dernière construction faisant intervenir la réalisation de Betti en plus de la réalisation de de Rham (ou $\omega$-réalisation), elle permet de prendre en compte le caractère multivalué de la fonction~$\Li_n$. On renvoie le lecteur à l'article de survol de \textcite{brownnotes} pour un tour d'horizon de ces notions.
		\end{rema}
		
		\begin{rema}\label{rema: relations periodes motiviques}
		\`{A} bien des égards, les relations entre polylogarithmes motiviques (et plus généralement entre périodes motiviques) se prouvent comme les relations entre polylogarithmes, c'est-à-dire par différentiation, le cocrochet jouant le rôle de la différentielle. Soit $B$ un ouvert non vide de $\mathbb{A}^1_F=\mathrm{Spec}(F[t])$. Par localisation, la $K$-théorie de $B$, définie comme la $K$-théorie de l'anneau des fonctions $\mathcal{O}(B)$, s'exprime en termes de la $K$-théorie de $F$ et de ses extensions finies. Par le théorème de Borel, on dispose donc si $F$ est un corps de nombres d'une catégorie $\MT(B)$ de motifs de Tate mixtes (ou faisceaux motiviques de Tate mixtes) sur $B$, et on montre que le morphisme structurel~$B\to \mathrm{Spec}(F)$ induit un isomorphisme 
		\begin{equation}\label{eq: rigidite}
		\mathrm{Ext}^1_{\MT(F)}(\QQ(-n),\QQ(0))\simeq \mathrm{Ext}^1_{\MT(B)}(\QQ(-n),\QQ(0))
		\end{equation}
		pour tout $n\geq 2$ (phénomène de \emph{rigidité}). Soit $\xi(t)\in \QQ[\mathcal{O}(B)^\times]$ tel que $\delta(\Li_n^{\mathcal{C}}(\xi(t)))=0$, avec $n\geq 2$. Alors par~\eqref{eq: rigidite} on voit que $\Li_n^{\mathcal{C}}(\xi(t))$ est indépendant de $t$, et il suffit de spécialiser $t$ en un point rationnel de $B$ pour trouver la constante. Si $F=\QQ$ et $n$ est pair, l'annulation de $K_{2n-1}(\QQ)_\QQ$ implique que la constante est nécessairement zéro. En travaillant avec des variations de structures de Hodge mixtes on peut traiter le cas de bases plus générales (voir par exemple \textcite{goncharovperiodsmixedmotives}).
		\end{rema}
		
		\begin{rema}
		Dans l'esprit de la remarque précédente, on peut voir $\mathrm{Li}_n$ comme un motif de Tate mixte  sur la base $\mathbb{P}^1\setminus \{0,1,\infty\}$. En réalisation de Hodge, on obtient une variation de structures de Hodge mixtes sur $\mathbb{P}^1(\CC)\setminus \{0,1,\infty\}$ étudiée par \textcite{hainclassical} (le système local sous-jacent avait été calculé par \textcite{ramakrishnan}). La présence d'une connexion (de Gauss--Manin) dans ce formalisme explique la ressemblance formelle entre l'équation différentielle~\eqref{eq: equa diff Li} et l'expression du cocrochet~\eqref{eq: cocrochet Li n}.
		\end{rema}
		
	\subsection{Interprétation motivique de la conjecture de Zagier}
		
		Le raisonnement de la remarque~\ref{rema: relations periodes motiviques} implique que le morphisme
		$$\iota_n:\mathcal{B}_n(F)\to \mathcal{C}_n(F) \; , \; [x]_n \mapsto \mathrm{Li}^{\mathcal{C}}_n(x)$$
		est bien défini pour tout $n\geq 2$. De plus, la formule~\eqref{eq: cocrochet Li n} implique que $\log^{\mathcal{C}}:F^\times_\QQ\stackrel{\sim}{\to} \mathcal{C}_1(F)$ et $\iota_2,\ldots,\iota_n$ s'assemblent en un morphisme de complexes
		$$\iota^\bullet(n): \mathcal{B}^\bullet(F,n) \To \mathrm{CE}^\bullet(\mathcal{C}(F))_n \ .$$
		On obtient donc, au niveau des $H^1$, un morphisme
		$$\psi_n : \ker(\delta_n) \To \ker\big(\delta:\mathcal{C}_n(F)\to (\Lambda^2\mathcal{C}(F))_n\big)\ .$$ 
		Si $F$ est un corps de nombres, on a donc grâce à~\eqref{eq: isom K theorie primitifs C} et~\eqref{eq: periode reelle Li n} le diagramme commutatif suivant :
		\begin{equation}\label{eq: diag comm zagier k theorie bis}
		\xymatrixcolsep{5pc}\xymatrix{
		\ker(\delta_n) \ar[r]^-{\psi_n} \ar@/_2pc/[rr]_-{\mathrm{P}_n^F} & K_{2n-1}(F)_\QQ \ar[r]^-{\varpi_n} & \RR^{d_n}
		}
		\end{equation}
		La ressemblance formelle avec le diagramme~\eqref{eq: diag comm zagier k theorie} ne doit pas tromper le lecteur sur le fait que $\psi_n$ va dans le sens opposé à $\varphi_n$ ! L'existence du diagramme~\eqref{eq: diag comm zagier k theorie bis} n'est pas suffisante pour prouver la conjecture de Zagier, qui serait une conséquence de la conjecture suivante (et de l'égalité des régulateurs $\varpi_n=\rho_n$).
		
		\begin{conj}\label{conj: zagier k theorie bis}
		Le morphisme $\psi_n$ est surjectif.
		\end{conj}
		
		En effet, \eqref{eq: zagier conj} découle alors de~\eqref{eq: borel theorem} et~\eqref{eq: diag comm zagier k theorie bis} en prenant pour $(\xi_j)_{1\leq j\leq d_n}$ des antécédents par $\psi_n$ d'une base de $K_{2n-1}(F)_\QQ$.
		
		\begin{rema}
		Le diagramme~\eqref{eq: diag comm zagier k theorie bis} a été construit par \textcite{beilinsondeligne} en supposant l'existence du formalisme motivique, et indépendamment par \textcite{dejeu}\footnote{\textcite[Remark 5.4]{dejeu} construit en fait une version \emph{entière} du morphisme $\psi_n$, ce qui montre que l'image de (la version entière de) $\mathrm{ker}(\delta_n)$ par $\mathrm{P}^F_n$ est un groupe abélien de type fini, comme le prévoit la conjecture originale de \textcite{zagierconjecture}.}. En fait, Beilinson--Deligne comme de Jeu rendent $\iota_n$ et donc $\psi_n$ injectifs en remplaçant l'espace $\mathcal{R}_n(F)$ des relations polylogarithmiques par l'espace des relations satisfaites par les $\Li_n^{\mathcal{C}}(x)$.
		\end{rema}
		
		\begin{rema}
		Via le dictionnaire tannakien, la conjecture~\ref{conj: zagier k theorie bis} affirme que toutes les extensions de~$\QQ(-n)$ par~$\QQ(0)$ dans~$\MT(F)$ se trouvent dans la sous-catégorie tannakienne engendrée par les motifs polylogarithmiques~$L(x)$ pour~\mbox{$x\in F^\times$}.
		\end{rema}
	
		Il est important de noter que le morphisme $\iota_n:\mathcal{B}_n(F)\to \mathcal{C}_n(F)$ lui-même n'est pas surjectif pour $n\geq 4$ (c'est-à-dire que les motifs polylogarithmiques n'engendrent pas la catégorie tannakienne $\MT(F)$). En effet, on peut par exemple produire des éléments de $\mathcal{C}_4(F)$ dont le cocrochet a une composante non nulle dans $\mathcal{C}_2(F)\wedge \mathcal{C}_2(F)$. C'est ce phénomène qui rend le cas $n=4$ de la conjecture de Zagier si crucial dans la perspective du cas général.
	
	\subsection{La conjecture de liberté}
	
		\subsubsection{\'{E}noncé}
	
		La conjecture suivante de \textcite{goncharovgeometry, goncharovseattle} renforce la conjecture~\ref{conj: zagier k theorie bis} de manière cohérente avec la conjecture~\ref{conj: goncharov polylog complexes}.
		
		\begin{conj}\label{conj: qis B C}
		Pour tout $n\geq 1$, le morphisme de complexes $\iota^\bullet(n)$ est un quasi-isomorphisme.
		\end{conj}
		
		On l'appelle \emph{conjecture de liberté} à cause de la formulation équivalente suivante.
	
		\begin{conj}\label{conj: liberte}
		\begin{enumerate}[(i)]
		\item La cogèbre de Lie quotient $\mathcal{C}(F)/\mathcal{C}_1(F)$ est colibre. 
		\item Pour tout $n\geq 2$, le morphisme $\iota_n$ est injectif et identifie $\mathcal{B}_n(F)$ avec l'espace des primitifs de $\mathcal{C}(F)/\mathcal{C}_1(F)$ en poids $n$.
		\end{enumerate}
		\end{conj}
		
		Dualement, cette conjecture affirme que la sous-algèbre de Lie $\mathfrak{u}(F)_{\leq -2}$ de l'algèbre de Lie motivique~$\mathfrak{u}(F)$ est libre et que ses indécomposables en poids $-n\leq -2$ sont donnés par l'espace $\mathcal{B}_n(F)^\vee$. Si $F$~est un corps de nombres, on a vu que l'algèbre de Lie motivique~$\mathfrak{u}(F)$ est libre et il en est donc de même pour $\mathfrak{u}(F)_{\leq -2}$, ce qui règle le point~\emph{(i)} .
		
		\begin{rema}
		La conjecture de liberté peut être séparée en deux énoncés relativement indépendants. D'une part, l'injectivité de $\iota_n$ revient à dire que les seules relations satisfaites par les polylogarithmes multiples motiviques $\Li^{\mathcal{C}}_n(x)$ sont celles prévues par l'espace des relations polylogarithmiques $\mathcal{R}_n(F)$. D'autre part, l'enjeu est de préciser la place des polylogarithmes motiviques à l'intérieur de la cogèbre de Lie motivique ; la conjecture implique notamment qu'on peut identifier via $\iota_n$ :
		\begin{equation}\label{eq: B inside C}
		\mathcal{B}_n(F) \stackrel{?}{=} \{ \alpha \in \mathcal{C}_n(F) \; |\; \delta \alpha \in \mathcal{C}_{n-1}(F)\wedge \mathcal{C}_1(F) \}\ .
		\end{equation}
		\end{rema}
		
		\subsubsection{La conjecture de liberté en bas poids}\label{par: liberte bas poids}
		
		Explicitons la conjecture de liberté en poids $n=2,3,4$.
		\begin{enumerate}[---]
		\item En poids $2$ elle affirme que $\iota_2:\mathcal{B}_2(F)\to\mathcal{C}_2(F)$ est un isomorphisme, ce qui découle des travaux de \textcite{suslinbloch} (voir la remarque \ref{rema: suslin bloch}).
		\item En poids $3$ elle affirme que $\iota_3:\mathcal{B}_3(F)\to\mathcal{C}_3(F)$ est un isomorphisme ; même dans le cas d'un corps de nombres, ni l'injectivité ni la surjectivité ne semblent être connues.
		\item En poids $4$ elle prend la forme d'une suite exacte courte
		\begin{equation}\label{eq: suite exacte courte liberte 4}
		0\To \mathcal{B}_4(F) \stackrel{\iota_4}{\To} \mathcal{C}_4(F) \stackrel{\delta_{2,2}}{\To} \mathcal{C}_2(F)\wedge \mathcal{C}_2(F)\To 0 \ .
		\end{equation}
		On sait seulement montrer que $\delta_{2,2}$ est surjectif si $F$ est un corps de nombres (voir \S\ref{par: relation 5 termes poids 4} plus bas).
		\end{enumerate}
	
	\subsection{Une stratégie vers la conjecture de Zagier, version raffinée}
		
		Le diagramme~\eqref{eq: diag comm zagier k theorie bis} semble être l'\emph{habitat naturel} de la conjecture de Zagier, et la conjecture~\ref{conj: zagier k theorie bis} (surjectivité de $\psi_n$) semble être plus fondamentale que la conjecture~\ref{conj: zagier k theorie}, qui nécessite d'inventer un morphisme $\varphi_n$ (conjecturalement une section de $\psi_n$) en découvrant des relations polylogarithmiques --- la conjecture~\ref{conj: zagier k theorie bis} n'étant pas concernée par ces relations.  Malheureusement, cette approche \og naturelle \fg{} à la conjecture de Zagier n'a pas encore porté ses fruits. L'obstacle principal est que la définition des motifs de Tate mixtes est si peu explicite qu'on ne sait pas très bien dire \og à quoi ressemble \fg{} un élément de $\mathcal{C}_n(F)$.
		
		Malgré ce constat, le point de vue motivique permet de raffiner la stratégie vers la conjecture de Zagier en mettant en lumière le rôle central joué par la cogèbre de Lie motivique.  Il semble en effet plus raisonnable de chercher à définir un morphisme de la~$K$-théorie vers (les éléments primitifs de) $\mathcal{C}(F)$, et de s'en remettre à la conjecture de liberté pour en déduire une flèche vers $\mathcal{B}(F)$ comme dans la conjecture~\ref{conj: zagier k theorie}.
		
		La conjecture de liberté étant présentement hors de portée, on procède via des remplacements \og symboliques \fg{} des espaces $\mathcal{B}_n$ et $\mathcal{C}_n$, qu'on note $B_n$ et $C_n$, et qui sont définis par générateurs et relations \emph{explicites}, avec la même structure abstraite (cocrochets) que leurs versions calligraphiques et des morphismes $B\to \mathcal{B}$ et $C\to \mathcal{C}$ (dont on peut conjecturer qu'ils sont des isomorphismes). Il faut penser à $B_n(F)$ comme à un remplacement de $\mathcal{B}_n(F)$ où les relations polylogarithmiques inexplicites $\mathcal{R}_n(F)$ sont remplacées par des familles explicites de relations. Le groupe $C_n(F)$ joue un rôle similaire mais incorpore la structure motivique de tous les polylogarithmes \emph{multiples}, qu'on étudiera au prochain paragraphe. La structure conjecturale de la cogèbre de Lie motivique joue un rôle de guide et suggère une stratégie vers le conjecture de Zagier découpée en deux étapes relativement indépendantes :
		
		\begin{enumerate}[1)]
		\item La première étape consiste à construire un morphisme 
		$$K_{2n-1}(F)_\QQ\To H^1(\mathrm{CE}^\bullet(C(F))_n)=\ker\big(\delta:C_n(F)\to (\Lambda^2C(F))_n\big)$$
		à partir d'un morphisme de complexes 
		\begin{equation}\label{eq: morphisme BG to CEC}
		BG^\bullet(F,n)\To \mathrm{CE}^\bullet(C(F))_n
		\end{equation}
		comme au \S\ref{par: strategie premiere}. Les techniques générales développées par Goncharov permettent de démontrer la compatibilité souhaitée entre le régulateur de Borel et le régulateur de Hodge induit sur $C_n(F)$ via le morphisme vers $\mathcal{C}_n(F)$. On obtient comme corollaire de cette première étape une version \og faible \fg{} de la conjecture de Zagier, où les polylogarithmes classiques sont remplacés par les polylogarithmes multiples. \textcite{goncharovarakelov} a apporté une contribution importante dans cette direction en décrivant le régulateur de Borel via une fonction appelée polylogarithme grassmannien (univalué). \textcite{charltonganglradchenkograssmannian} donnent une formule pour une variante multivaluée de cette fonction, introduite par \textcite{goncharovsimple}, en termes de polylogarithmes multiples. En principe, la marche à suivre est donc d'étendre ces formules en un morphisme de complexes \eqref{eq: morphisme BG to CEC}. Nous ne rentrerons pas dans les détails de cette partie du travail de \textcite{goncharovrudenko}, dans le cas $n=4$, qui s'appuie notamment sur des travaux antérieurs de \textcite{goncharovgeometrytrilog}.
		\item La deuxième étape consiste à prouver une version \og symbolique \fg{} de la conjecture de liberté, c'est-à-dire à démontrer que le morphisme naturel $B^\bullet(F,n)\To \mathrm{CE}^\bullet(F,n)$ est un quasi-isomorphisme\footnote{Pour l'application à la conjecture de Zagier, on a seulement besoin de montrer qu'il induit un isomorphisme au niveau des $H^1$.}.  Pour $n=4$, cette tâche a été accomplie par \textcite{ganglfour} et par \textcite{goncharovrudenko}, et nous en décrirons les grandes lignes après un interlude sur les polylogarithmes multiples dans le contexte motivique.
		\end{enumerate}

\section{Polylogarithmes multiples motiviques et profondeur}
		
		On explique maintenant comment définir des versions motiviques des polylogarithmes multiples dans la cogèbre de Lie motivique, 
		$$\Li_{n_1,\ldots,n_r}^{\mathcal{C}}(x_1,\ldots,x_r) \in \mathcal{C}_{n_1+\cdots+n_r}(F) \ ,$$
		pour des entiers $n_1,\ldots,n_r\geq 1$ et des éléments $x_1,\ldots,x_r\in F^\times$. Ils sont définis comme des coefficients matriciels d'un (ind-)motif de Tate mixte sur $F$ appelé \emph{groupe fondamental motivique}, ce qui reflète l'interprétation des polylogarithmes multiples comme intégrales itérées. Le concept de profondeur (l'entier $r$) fera l'objet d'une discussion particulière. On introduit aussi certaines variantes, les \emph{corrélateurs}, en suivant \textcite{goncharovcorrelators}.
		
		\subsection{Intégrales itérées}
		
			Pour une variété différentielle $X$, des $1$-formes différentielles $\omega_1,\ldots,\omega_n$ sur $X$, et un chemin lisse $\gamma:[0,1]\to X$, on définit l'\emph{intégrale itérée}
			$$\int_\gamma \omega_1\cdots \omega_n = \int_{0\leq t_1\leq\cdots \leq t_n\leq 1} f_1(t_1)dt_1\cdots f_n(t_n)dt_n$$
			où l'on a noté $f_i(t)dt$ le tiré en arrière de $\omega_i$ par $\gamma$. Cela généralise le concept classique d'intégrale d'une $1$-forme le long d'un chemin à des mots formés de $1$-formes.
			
			Dans le cas où $X=\mathbb{P}^1(\CC)\setminus \{a_1,\ldots,a_n,\infty\}$ est la droite projective épointée, on note, pour des points-base $a_0,a_{n+1}\in X$, 
			\begin{equation}\label{eq: integrale iteree}
			\mathrm{I}(a_0;a_1,\ldots,a_n;a_{n+1}) = \int_{a_0}^{a_{n+1}} \frac{dt}{t-a_1}\cdots \frac{dt}{t-a_n}
			\end{equation}
			qui dépend d'un choix de chemin de $a_0$ vers $a_{n+1}$ dans $X$. On peut étendre cette définition au cas où $a_0, a_{n+1}$ sont \og à l'infini \fg{} dans $X$, c'est-à-dire parmi $a_1,\ldots,a_n,\infty$. Dans ce cas-là l'intégrale itérée peut diverger (ce qui est le cas si et seulement si $a_0=a_1$ ou $a_n=a_{n+1}$) et on fixe des vecteurs tangents à $\mathbb{P}^1(\CC)$ en $a_0$ et $a_{n+1}$ (qui sont alors appelés \emph{points-base tangentiels}) pour la régulariser. Il est alors naturel de considérer des chemins $\gamma:[0,1]\to\mathbb{P}^1(\CC)$ qui partent de $a_0$ et arrivent en $a_{n+1}$ avec des vecteurs vitesse égaux aux vecteurs tangents prescrits, et sont tels que $\gamma(t)\in X$ pour $0<t<1$.
			
			Pour des entiers $n_1,\ldots,n_r\geq 1$ et des nombres complexes non nuls $a_1,\ldots,a_r$ on note
			$$\mathrm{I}_{n_1,\ldots,n_r}(a_1,\ldots,a_r;a_{r+1}) = (-1)^r\, \mathrm{I}(0;a_1,\{0\}^{n_1-1},a_2,\{0\}^{n_2-1},\ldots,a_r,\{0\}^{n_r-1};a_{r+1})$$
			où $\{0\}^k$ désigne une liste de $0$ de longueur $k$. Cette notation est abusive et dépend du choix d'un chemin de $0$ à $a_{r+1}$. Dans le cas divergent $(n_r=1, a_r=a_{r+1})$ la convention est de régulariser l'intégrale itérée via le vecteur tangent $a_{r+1}$ en $a_{r+1}$, de sorte que $\operatorname{I}_1(a;a)=0$. 
			
			Un calcul rapide montre que les polylogarithmes multiples s'expriment en termes d'intégrales itérées sous la forme :
			$$\Li_{n_1,\ldots,n_r}(x_1,\ldots,x_r) = \mathrm{I}_{n_1,\ldots,n_r}( 1,x_1,x_1x_2,\ldots,x_1x_2\cdots x_{r-1} ; x_1x_2\cdots x_r)\ ,$$
			ou réciproquement :
			$$\mathrm{I}_{n_1,\ldots,n_r}(a_1,\ldots,a_r;a_{r+1}) = \Li_{n_1,\ldots,n_r}(\textstyle\frac{a_2}{a_1},\frac{a_3}{a_2},\ldots,\frac{a_{r+1}}{a_r})\ .$$
						
		\subsection{Groupoïdes fondamentaux motiviques et polylogarithmes multiples motiviques}
		
			Soit $F$ un corps de nombres pour lequel on fixe un plongement complexe. Soit $S\subset \mathbb{P}^1(F)$ un ensemble fini contenant $\infty$, et notons $X=\mathbb{P}^1_F\setminus S$. On fixe un ensemble fini de points-base (éventuellement tangentiels) sur $X$, définis sur $F$. On note~$\pi_1^{\mathrm{uni}}(X(\CC))$ la \emph{complétion pro-unipotente} (ou \emph{de Malcev}) du groupoïde fondamental topologique de $X(\CC)$ relatif au choix de points-base, qui est un schéma en groupoïdes défini sur $\mathbb{Q}$ qui factorise les représentations unipotentes de $\pi_1(X(\CC))$. Les intégrales itérées donnent lieu, d'après \textcite{cheniterated}, à un isomorphisme :
			\begin{equation}\label{eq: chen}
			T(H^1_{\mathrm{dR}}(X))\otimes_F\CC \stackrel{\sim}{\To} \mathcal{O}(\pi_1^{\mathrm{uni}}(X(\CC))_{a,b})\otimes_\QQ\CC\ ,
			\end{equation}
			où $T(V)$ désigne l'algèbre tensorielle sur $V$ et la flèche associe à un mot $\omega_1\cdots \omega_n$ la fonction d'intégration $\gamma\mapsto \int_\gamma\omega_1\cdots \omega_n$ définie sur les chemins $\gamma$ de $a$ à $b$ dans $X(\CC)$. On peut, d'après \textcite{deligneP1}, \textcite{goncharovdihedral, goncharovgaloissymmetries} et \textcite{delignegoncharov}, voir~\eqref{eq: chen} comme l'isomorphisme des périodes d'un ind-objet de la catégorie $\MT(F)$, noté $\mathcal{O}(\pi_1^{\mathrm{mot}}(X)_{a,b})$. En réalisation de Hodge cette construction remonte aux travaux de \textcite{morgantopology} et de \textcite{hainfundamentalgroup}.
			
			En $\omega$-réalisation, le $\QQ$-espace vectoriel $\mathrm{Hom}_{\MT(F)}(\QQ(-k),\mathrm{gr}_{2k}^W \mathcal{O}(\pi_1^{\mathrm{mot}}(X)_{a,b}))$ est indépendant de $a,b$ et a une base formée des mots de longueur $k$ en les $dt/(t-s)$, pour~$s\in S\setminus \{\infty\}$. On a donc un isomorphisme canonique $\mathrm{gr}_0^W\mathcal{O}(\pi_1^{\mathrm{mot}}(X)_{a,b}) \to \QQ(0)$, et l'intégrande de~\eqref{eq: integrale iteree} permet de définir un coefficient matriciel
			$$\mathrm{I}^{\mathcal{C}}(a_0;a_1,\ldots,a_n;a_{n+1}) \in \mathcal{C}_n(F)\ ,$$
			pour tous les $a_0,\ldots,a_{n+1}\in F$, qui ne dépend pas d'un choix de chemin de $a_0$ à $a_{n+1}$, et ne dépend de choix de vecteurs tangents en $a_0$ et $a_{n+1}$ que si $n=1$. 
			Par les formules du paragraphe précédent, on obtient aussi des éléments		
			\begin{equation}\label{eq: polylog multiples motiviques}
			\mathrm{I}_{n_1,\ldots,n_r}^{\mathcal{C}}(a_1,\ldots,a_r;a_{r+1}) \in \mathcal{C}_{n_1+\cdots+n_r}(F) \;\hbox{ et }\; \Li^{\mathcal{C}}_{n_1,\ldots,n_r}(x_1,\ldots,x_r)\in \mathcal{C}_{n_1+\cdots+n_r}(F)\ ,
			\end{equation}
			associés à $a_1,\ldots,a_{r+1}$ ou $x_1,\ldots,x_r$ dans $F^\times$. Ils généralisent~\eqref{eq: polylog motivique} et on les appelle \emph{polylogarithmes multiples motiviques} (dans leur version modulo produits). La conjecture suivante est due à \textcite{goncharovICM}.		
		
			\begin{conj}\label{conj: universalite}
			Les polylogarithmes multiples motiviques~\eqref{eq: polylog multiples motiviques} engendrent la cogèbre de Lie motivique $\mathcal{C}(F)$.
			\end{conj}
			
			Via le dictionnaire tannakien, cette conjecture affirme que la catégorie tannakienne~$\MT(F)$ est engendrée par les ind-objets $\mathcal{O}(\pi_1^{\mathrm{mot}}(\mathbb{P}^1_F\setminus S)_{a,b})$, pour $S\subset \mathbb{P}^1(F)$ un sous-ensemble fini contenant~$\infty$ et~$a,b$ des points-base (éventuellement tangentiels) définis sur~$F$.

			\begin{rema}
			Parmi les progrès vers cette conjecture, citons le théorème de  \textcite{brownMTZ} selon lequel l'objet $\mathcal{O}(\pi_1^{\mathrm{mot}}(\mathbb{P}^1_\QQ\setminus \{0,1,\infty\}))$, avec points-base tangentiels adéquats, engendre la sous-catégorie $\MT(\ZZ)\subset \MT(\QQ)$ des motifs de Tate mixtes sur $\mathbb{Z}$, et les résultats antérieurs de~\textcite{delignecyclotomic} dans le cas de certains corps cyclotomiques. Dans une autre direction,  \textcite{boehm} montre que tous les volumes de polytopes hyperboliques s'expriment en termes de polylogarithmes multiples (voir \textcite{rudenkodepth}). 
			\end{rema}
			
			\textcite{goncharovmultiple, goncharovgaloissymmetries} montre que le cocrochet des polylogarithmes multiples motiviques est calculé par la jolie formule suivante\footnote{Le lecteur attentif remarquera que certains termes du membre de droite peuvent dépendre de choix de vecteurs tangents en $a_1,\ldots,a_n$, même si la somme n'en dépend pas.}, qui généralise~\eqref{eq: cocrochet Li n} :
			\begin{equation}\label{eq: cocrochet integrale iteree}
			\begin{split}
			\delta \operatorname{I}^{\mathcal{C}}(a_0;&a_1,\ldots,a_n;a_{n+1})  \\
			&= \sum_{0\leq i<j\leq n} \operatorname{I}^{\mathcal{C}}(a_0;a_1,\ldots,a_i,a_{j+1},\ldots,a_n;a_{n+1})
			\wedge \operatorname{I}^{\mathcal{C}}(a_i;a_{i+1},\ldots,a_{j-1};a_j)\ .
			\end{split}
			\end{equation}

	\subsection{La conjecture de profondeur}
		
		On définit de manière récursive une filtration croissante $D$ (par la \emph{profondeur}) sur la cogèbre de Lie motivique, en posant $D_{-1}\mathcal{C}(F)=0$, $D_0\mathcal{C}(F)=\mathcal{C}_1(F)$ et pour $k\geq 1$,
		$$D_k\mathcal{C}(F) = \{\alpha \in \mathcal{C}(F) \; | \; \delta\alpha \in \mathcal{C}_1(F)\wedge \mathcal{C}(F) + D_{k-1}\mathcal{C}(F)\wedge D_{k-1}\mathcal{C}(F)\}\ .$$
		Il s'agit d'une filtration exhaustive : pour tout entier $n$ on a 
		\begin{equation}\label{eq: profondeur moitie}
		\mathcal{C}_n(F) = D_{\lfloor n/2 \rfloor}\mathcal{C}_n(F)\ .
		\end{equation}
		Cette définition est plus claire dans le contexte dual de l'algèbre de Lie motivique : la filtration (décroissante) duale est donnée par $D^0\mathfrak{u}(F)=\mathfrak{u}(F)$, $D^1\mathfrak{u}(F)=\mathfrak{u}_{\leq -2}(F)$, et pour tout $k\geq 2$, $D^k\mathfrak{u}(F)=[\mathfrak{u}_{\leq -2}(F),D^{k-1}\mathfrak{u}(F)]$.
		
		La filtration par la profondeur est compatible avec la notion de profondeur des polylogarithmes multiples, au sens où on a $\Li_{n_1,\ldots,n_r}^{\mathcal{C}}(x_1,\ldots,x_r)\in D_r\mathcal{C}(F)$ pour tous les~$x_1,\ldots,x_r\in F^\times$, comme on peut le voir en utilisant~\eqref{eq: cocrochet integrale iteree}. La conjecture suivante de \textcite{goncharovmultiple}, appelée \emph{conjecture de profondeur}, raffine la conjecture~\ref{conj: universalite} et précise la place des polylogarithmes multiples motiviques dans la cogèbre de Lie motivique.
		
		\begin{conj}\label{conj: profondeur}
		Soit $k\geq 1$ un entier. Les polylogarithmes multiples motiviques $\Li_{n_1,\ldots,n_r}^{\mathcal{C}}(x_1,\ldots,x_r)$, pour $r\leq k$, engendrent $D_k\mathcal{C}(F)$.
		\end{conj}
		
		Pour $k=1$ on retrouve~\eqref{eq: B inside C}. Le théorème suivant a été prouvé récemment par \textcite{rudenkodepth} et constitue, au vu de~\eqref{eq: profondeur moitie}, un premier pas vers la conjecture de profondeur.
		
		\begin{theo}\label{theo: rudenko depth}
		Tout polylogarithme multiple motivique de poids $n\geq 2$ peut s'écrire comme combinaison linéaire de polylogarithmes multiples motiviques $\Li_{n_1,\ldots,n_r}^{\mathcal{C}}(x_1,\ldots,x_r)$ avec $r\leq \lfloor n/2 \rfloor$.
		\end{theo}
		
		On peut enlever le mot \og motivique \fg{} de cet énoncé quitte à travailler modulo des produits de polylogarithmes multiples de poids inférieur. Ce théorème est une amélioration notable de l'état antérieur des connaissances, où la borne $\lfloor n/2\rfloor$ était remplacée par $\max\{1,n-2\}$.
		
	\subsection{Corrélateurs}
			
		\textcite{goncharovcorrelators} a introduit une manière d'organiser les polylogarithmes multiples motiviques sous une forme plus symétrique (au sens de la symétrie cyclique voire diédrale) qui semble bien adaptée à l'étude des relations. Les objets centraux sont appelés \emph{corrélateurs motiviques}.
			
		\subsubsection{Produit et crochet d'Ihara}
				
		Soit $S\subset \mathbb{P}^1(F)$ un ensemble fini contenant $\infty$, et notons $S\setminus \{\infty\}=\{s_1,\ldots,s_N\}$. On fait le choix d'un vecteur tangent à~$\mathbb{P}^1_F$ en chaque point de~$S$ et on considère le groupoïde fondamental motivique~$\pi_1^{\mathrm{mot}}(\mathbb{P}^1_F\setminus S)$ relatif à ce choix de points-base. On note~$\Pi$ sa $\omega$-réalisation, qui est un schéma en groupoïdes défini sur~$\QQ$. Par le formalisme tannakien, il est muni d'une action du groupe tannakien~$U(F)$, qui a d'abord été étudiée par \textcite{iharagalois} dans le cadre des représentations galoisiennes. Notons d'abord que le groupoïde~$\Pi$ est constant : pour une $\QQ$-algèbre~$R$ et~$a,b\in S$, $\Pi_{a,b}(R)$ est canoniquement isomorphe au groupe des séries non commutatives \og group-like \fg{} en~$N$ variables $X_1,\ldots,X_N$ :
		$$\Pi_{a,b}(R) = \{ F \in R\langle\langle X_1,\ldots,X_N\rangle\rangle \; ,\; \Delta(F)=F\otimes F\}\ ,$$
		où le coproduit (complété) $\Delta$ est déterminé par $\Delta(1)=1\otimes 1$ et $\Delta(X_i)=1\otimes X_i+X_i\otimes 1$. Cependant, l'action de $U(F)$ sur $\Pi_{a,b}$ est sensible aux points-base $a,b$. 
		
		Soit $A$ le schéma en groupes dont les points sont les automorphismes du groupoïde $\Pi$ qui agissent trivialement sur $\exp(X_i)\in \Pi_{s_i,s_i}$ pour tout $i=1,\ldots,N$ et sur $\exp(X_\infty)\in \Pi_{\infty,\infty}$, avec $X_\infty=-X_1-\cdots-X_N$. L'action de $U(F)$ sur $\Pi$ se factorise par $A$ et induit donc un morphisme de schémas en groupes $U(F)\To A$. L'évaluation en la série constante $1\in \Pi_{\infty,s_i}$, pour $i=1,\ldots,N$, induit un isomorphisme de schémas $A\stackrel{\sim}{\To} P$ où $P$ est le sous-schéma de $\Pi_{\infty,s_1}\times \cdots\times \Pi_{\infty,s_N}$ formé des $N$-uplets $(F_1,\ldots,F_N)$ qui vérifient $F_1X_1F_1^{-1}+\cdots +F_NX_NF_N^{-1} = X_1+\cdots +X_N$. On a donc par transfert une loi de groupe sur $P$ qu'on appelle \emph{produit d'Ihara}. 
		
		Passons maintenant aux algèbres de Lie en notant $\mathfrak{a}=\mathrm{Lie}(A)$ et $\mathfrak{p}=\mathrm{Lie}(P)$. Cette dernière algèbre de Lie est l'espace des $N$-uplets $(f_1,\ldots,f_N)$ d'éléments de l'algèbre de Lie libre complétée $\widehat{\mathrm{Lie}}(X_1,\ldots,X_N)$ qui vérifient $[f_1,X_1]+\cdots + [f_N,X_N]=0$. Son crochet de Lie, noté $\{-,-\}$ et appelé \emph{crochet d'Ihara}, est donné par la formule $\{\underline{f},\underline{g}\}=\underline{h}$ avec
		$$h_i=D_{\underline{f}}(g_i)-D_{\underline{g}}(f_i)-[f_i,g_i]\ ,$$
		où l'on note $D_{\underline{f}}$ l'unique dérivation continue de $\widehat{\mathrm{Lie}}(X_1,\ldots,X_N)$ qui est telle que $D_{\underline{f}}(X_i)=[f_i,X_i]$ pour tout $i=1,\ldots,N$.
		
		\subsubsection{Codage par des mots cycliques}
		
		\textcite{drinfeld} donne la description suivante de l'algèbre de Lie $\mathfrak{p}$. On voit $\mathfrak{p}$ à l'intérieur de l'algèbre de Lie $\mathfrak{q}$, qui est définie de la même manière en permettant aux $f_i$ d'être des séries non commutatives quelconques dans $\QQ\langle\langle X_1,\ldots,X_N\rangle\rangle$. Soit maintenant $\mathfrak{cyc}$ l'espace des séries cycliques, c'est-à-dire invariantes par permutation circulaire des lettres à l'intérieur de chaque mot, qui ont un terme constant nul. Pour une série $F$ on note $\partial_i F$ sa dérivée partielle par rapport à $X_i$, qui est donnée sur les mots par $\partial_i(X_iw)=w$ et $\partial_i(X_jw)=0$ pour $j\neq i$. On a un isomorphisme $\mathfrak{cyc}\stackrel{\sim}{\To} \mathfrak{q}$ donné par $F\mapsto (\partial_1F,\ldots,\partial_nF)$. On a donc le diagramme suivant d'algèbres de Lie :
		\begin{equation}\label{eq: correlateurs algebres lie}
		\xymatrixcolsep{4pc}\xymatrix{
		\mathfrak{u}(F) \ar[r]& \mathfrak{a} \ar[r]^-{\sim} & \mathfrak{p} \ar@{^(->}[r]& \mathfrak{q} &  \mathfrak{cyc} \ar[l]_-{\sim}
		}
		\end{equation}
		
		Pour décrire le crochet de Lie sur l'algèbre de Lie $\mathfrak{cyc}$ il est plus commode de la voir comme le dual linéaire d'une cogèbre de Lie graduée $\mathcal{C}yc$. Pour tout entier $n\geq 0$, une base de $\mathcal{C}yc_n$ est donnée par les mots de longueur $n+1$ en $X^1,\ldots,X^N$ considérés modulo permutations circulaires des lettres à l'intérieur de chaque mot. Le cocrochet sur $\mathcal{C}yc$ est donné par
		$$\delta(X^{i_0}\cdots X^{i_n}) = \sum_{0\leq j<k\leq n} (X^{i_j}\cdots X^{i_k}) \wedge (X^{i_k}\cdots X^{i_{j-1}})\ .$$
		On a donc un morphisme de cogèbres de Lie 
		\begin{equation}\label{eq: morphisme Cyc to C}
		\mathcal{C}yc \to \mathcal{C}(F)\ ,
		\end{equation}
		dual à la composition des morphismes dans \eqref{eq: correlateurs algebres lie}.
		On montre que ce morphisme ne dépend des choix des vecteurs tangents en les points-base $s_1,\ldots,s_N, \infty$ qu'en poids $1$. 
		
		\subsubsection{Corrélateurs motiviques}		
		
\begin{figure}[h]
		\begin{center}
		\begin{tikzpicture}

		\def \r {2.5} ;
		\draw (0,0) circle (\r) ;
		\foreach \k in {0,...,8} \draw (90-360/9*\k:\r) node[circle, inner sep = 1pt, fill = black] (\k)  {};
		\foreach \k in {0,...,8} \node at (90-360/9*\k:1.2*\r) {$x_\k$};
		\draw  (90-360/9*5/2:\r) node[rectangle, inner sep = 2pt, fill = blue] (mid) {};
		\draw[thick, blue] (mid) -- (6);
		
		\end{tikzpicture}
		\end{center}
		\caption{Pour $n=8$, le terme dans la formule de cocrochet \eqref{eq: cocrochet correlateur} correspondant à $(j,k)=(3,6)$ est $\mathrm{Cor}^{\mathcal{C}}(x_3,x_4,x_5,x_6)\wedge \mathrm{Cor}^{\mathcal{C}}(x_6,x_7,x_8,x_0, x_1,x_2)$.}\label{fig: cocrochet correlateurs}
		\end{figure}
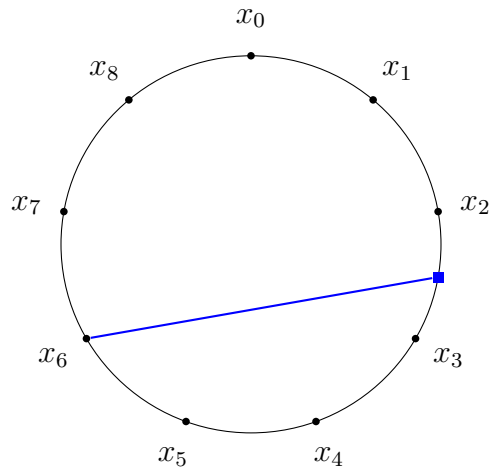

		Soient $x_0,x_1,\ldots,x_n\in F$ et appliquons la construction du paragraphe précédent à un ensemble $S\supset \{x_0,\ldots,x_n\}$. En appliquant \eqref{eq: morphisme Cyc to C} au mot dont les lettres correspondent à $x_0,\ldots,x_n$, on obtient un élément
		$$\mathrm{Cor}^{\mathcal{C}}(x_0,x_1,\ldots,x_n) \in \mathcal{C}_n(F)\ ,$$
		appelé \emph{corrélateur motivique} par \textcite{goncharovcorrelators}.
		On a la symétrie cyclique $\mathrm{Cor}^{\mathcal{C}}(x_0,x_1,\ldots,x_n)=\mathrm{Cor}^{\mathcal{C}}(x_1,\ldots,x_n,x_0)$ et la formule de cocrochet
		\begin{equation}\label{eq: cocrochet correlateur}
		\delta\operatorname{Cor}^{\mathcal{C}}(x_0,x_1,\ldots,x_n) = \sum_{j,k} \operatorname{Cor}^{\mathcal{C}}(x_j,\ldots,x_k) \wedge \operatorname{Cor}^{\mathcal{C}}(x_k,\ldots,x_{j-1})\ ,
		\end{equation}
		où la somme porte sur les indices $j,k$ tels que $k\neq j$ et $k\neq j-1$ modulo $n+1$ (ou plus visuellement, sur les coupures dans un disque, comme dans la figure \ref{fig: cocrochet correlateurs}).
				
		En poids $1$ on a 
		$$\mathrm{Cor}^{\mathcal{C}}(x_0,x_1) = \log^{\mathcal{C}}(x_1-x_0)$$
		si on fait le choix du vecteur tangent en $\infty$ qui est dual à $d(1/x)$. En poids $2$ on a
		$$\mathrm{Cor}^{\mathcal{C}}(x_0,x_1,x_2) = \Li_2^{\mathcal{C}}(\textstyle\frac{x_1-x_0}{x_2-x_0})\ .$$
		
		\begin{rema} 
		Goncharov donne une formule intégrale pour $\varpi_n(\operatorname{Cor}^{\mathcal{C}}(x_0,\ldots,x_n))$ et étudie des variantes sur des courbes de genre supérieur qui devraient jouer un rôle dans des généralisations de la conjecture de Zagier (voir la remarque \ref{rema: regulateurs}). 
		\end{rema}		
		
		Les corrélateurs motiviques s'expriment en termes des polylogarithmes multiples motiviques, et \emph{vice versa}. Par exemple, on montre facilement par récurrence sur $n\geq 1$ qu'on a l'égalité, pour $x\in F$ : 
		$$\Li_n^{\mathcal{C}}(x) = - \mathrm{Cor}^{\mathcal{C}}(\underbrace{0,\ldots,0}_{n-1},1,x)\ .$$
		Il semble néanmoins que les relations entre corrélateurs motiviques aient une combinatoire plus maniable que celles entre polylogarithmes multiples motiviques, notamment du fait de la symétrie cyclique (comparer les formules~\eqref{eq: cocrochet integrale iteree} et~\eqref{eq: cocrochet correlateur}). Ils sont à la source de la description \og symbolique \fg{} par \textcite{goncharovrudenko} de la cogèbre de Lie motivique en poids $\leq 4$, que nous abordons maintenant.

\section{La structure de la cogèbre de Lie motivique en poids $\leq 4$}

	On s'intéresse maintenant à la structure fine de la cogèbre de Lie motivique en poids~$\leq 4$ et aux résultats de \textcite{ganglfour} et de \textcite{goncharovrudenko} vers la conjecture de liberté et la conjecture de Zagier pour $n=4$.
		
	\subsection{La relation à $5$ termes en poids $4$}\label{par: relation 5 termes poids 4}
	
		Rappelons que la conjecture de liberté prévoit, en poids $4$, la suite exacte courte~\eqref{eq: suite exacte courte liberte 4}. En mettant de côté l'injectivité de $\iota_4$, cela revient à montrer que la composante $\delta_{2,2}$ du cocrochet induit un isomorphisme
		$$\overline{\delta}_{2,2} : \mathcal{C}_4(F)/\iota_4(\mathcal{B}_4(F)) \To \mathcal{C}_2(F)\wedge \mathcal{C}_2(F)\ .$$
		
		Rappelons (voir la remarque \ref{rema: suslin bloch} et le \S\ref{par: liberte bas poids}) que $\mathcal{C}_2(F)$ a une présentation par générateurs et relations, où les générateurs sont les dilogarithmes motiviques $\Li_2^{\mathcal{C}}(x)$ et les relations sont les relations à $5$ termes et leurs spécialisations. On obtient donc facilement la surjectivité de $\overline{\delta}_{2,2}$, puisque la formule générale \eqref{eq: cocrochet integrale iteree} implique l'égalité, pour $x,y\in F^\times$ :
		\begin{equation}\label{eq: delta 2 2 I 3 1}
		\delta_{2,2}\operatorname{I}^{\mathcal{C}}_{3,1}(x,y;1) = - \Li_2^{\mathcal{C}}(x)\wedge \Li_2^{\mathcal{C}}(y)\ .
		\end{equation}
		
		Pour prouver l'injectivité, il faut déjà montrer que tout polylogarithme multiple motivique de poids $4$ dont le cocrochet $\delta_{2,2}$ s'annule est une combinaison linéaire de polylogarithmes motiviques $\Li_4^{\mathcal{C}}(x)$. (C'est suffisant si l'on croit à la conjecture~\ref{conj: universalite}, ou dans une version \og symbolique \fg{} de la cogèbre de Lie motivique où seuls les polylogarithmes multiples apparaissent.) Une première réduction (voir \textcite{ganglfour} et \textcite{dan} corrigé par \textcite{charltondan}) permet d'exprimer tout polylogarithme multiple motivique de poids $4$ en termes de $\Li_4^{\mathcal{C}}$ et $\Li_{3,1}^{\mathcal{C}}$ seulement. Au vu de~\eqref{eq: delta 2 2 I 3 1} on est alors amené à définir une section de $\overline{\delta}_{2,2}$ par la formule $\Li_2^{\mathcal{C}}(x)\wedge \Li_2^{\mathcal{C}}(y)\mapsto -\operatorname{I}_{3,1}^{\mathcal{C}}(x,y;1)$, et le tout est de montrer que cette section est bien définie. On montre assez facilement que $\operatorname{I}_{3,1}^{\mathcal{C}}(x,y;1) + \operatorname{I}_{3,1}^{\mathcal{C}}(y,x;1)=0$ modulo $\iota_4(\mathcal{B}_4(F))$ et il reste à montrer que $\operatorname{I}^{\mathcal{C}}_{3,1}(x,y;1)$ vérifie la relation à~$5$ termes en la variable $x$. C'est le contenu du théorème suivant, conjecturé par \textcite{goncharovseattle, goncharovgeometry} et démontré par \textcite{ganglfour}.
		
		\begin{theo}\label{theo: gangl}
		On a, pour $x_0,\ldots,x_4\in \mathbb{P}^1(F)$ deux à deux distincts et $y\in F^\times$ :
		\begin{equation}\label{eq: five term gangl}
		\sum_{i=0}^4 (-1)^i \operatorname{I}^{\mathcal{C}}_{3,1}(r(x_0,\ldots,\widehat{x_i},\ldots,x_4),y;1) = 0 \;\;\; \mathrm{mod}\; \iota_4(\mathcal{B}_4(F))\ .
		\end{equation}
		\end{theo}
		
		Ce théorème apporte du poids à la conjecture de liberté (conjecture~\ref{conj: liberte}) et est une étape cruciale dans la preuve de la conjecture de Zagier pour $n=4$. La preuve de Gangl, assistée par ordinateur, exprime le membre de gauche de~\eqref{eq: five term gangl} comme une combinaison linéaire explicite de 122 évaluations de $\Li_4^{\mathcal{C}}$ en des arguments qui sont malheureusement difficiles à interpréter en vue d'éventuelles généralisations. L'approche de \textcite{goncharovrudenko} permet de redémontrer ce résultat de manière plus conceptuelle grâce à des relations entre polylogarithmes multiples motiviques qui semblent avoir une origine systématique.  Ces relations donnent lieu à une version \og symbolique \fg{} de la cogèbre de Lie motivique en poids $\leq 4$, qu'on présente maintenant.
		
	\subsection{Définitions de $B$ et $C$ en poids $\leq 4$}		
			
		En suivant \textcite{goncharovrudenko}, on définit maintenant par générateurs et relations des espaces vectoriels $C_n(F)$ et des sous-espaces $B_n(F)$ de manière fonctorielle en le corps $F$, pour $n=1,2,3,4$. On pose $B_1(F)=C_1(F)=F^\times_\QQ$.
		
		On considère la normalisation suivante du birapport de $4$ points dans $\mathbb{P}^1$
		$$r_2(x_1,x_2,x_3,x_4)=\frac{(x_1-x_2)(x_3-x_4)}{(x_2-x_3)(x_4-x_1)}\ ,$$
		ainsi qu'une version à $6$ points :
		$$r_3(x_1,x_2,x_3,x_4,x_5,x_6)=-\frac{(x_1-x_2)(x_3-x_4)(x_5-x_6)}{(x_2-x_3)(x_4-x_5)(x_6-x_1)} \ \cdot$$
		
		On voit $r_2$ (resp. $r_3$) comme un morphisme de $\overline{\mathcal{M}}_{0,4}$ (resp. $\overline{\mathcal{M}}_{0,6}$) dans $\mathbb{P}^1$. Rappelons que l'espace de modules $\mathcal{M}_{0,k}$ est l'espace des $k$-uplets de points deux à deux distincts de la droite projective $\mathbb{P}^1$ modulo les automorphismes de $\mathbb{P}^1$, et que $\overline{\mathcal{M}}_{0,k}$ désigne sa compactification de Deligne--Mumford.  Dans les définitions qui suivent on fait un abus de notation en notant $(x_1,\ldots,x_k)$ un point de $\overline{\mathcal{M}}_{0,k}$, et de manière similaire les points qui s'en déduisent par l'action du groupe symétrique sur $k$ éléments ou par des applications d'oubli $\overline{\mathcal{M}}_{0,k}\to \overline{\mathcal{M}}_{0,k'}$ avec $k'<k$.
	
		\begin{defi}
		On définit $B_2(F)=C_2(F)$ comme le quotient de l'espace vectoriel librement engendré par des symboles $[x]_2$ pour $x\in \mathbb{P}^1(F)$ par les relations suivantes, pour $(x_1,\ldots,x_5)\in \overline{\mathcal{M}}_{0,5}(F)$ :
		\begin{equation*}
		\begin{split} 
		[r_2(x_1,x_2,x_3,x_4)]_2  +  [r_2(x_2,x_3,&x_4,x_5)]_2 + [r_2(x_3,x_4,x_5,x_1)]_2  \\ & +[r_2(x_4,x_5,x_1,x_2)]_2 + [r_2(x_5,x_1,x_2,x_3))]_2= 0\ .
		\end{split}
		\end{equation*}
		\end{defi}
		
		On a ici réécrit la relation à $5$ termes~\eqref{eq: five term birapport} sous une forme cyclique. Dans la suite on utilise la sommation cyclique $\mathrm{Cyc}_nf(x_1,\ldots,x_n) = \sum_{i\in\mathbb{Z}/n\mathbb{Z}} f(x_i,x_{i+1},\ldots,x_{i+n})$, où les indices sont pris modulo $n$.
		
		\begin{defi}
		On définit $B_3(F)=C_3(F)$ comme le quotient de l'espace vectoriel engendré par des symboles $[x]_3$ et $[x,y]_{2,1}$ pour $x,y\in \mathbb{P}^1(F)$ par les relations suivantes, pour $(x_1,\ldots, x_6)\in \overline{\mathcal{M}}_{0,6}(F)$ :
		\begin{equation}\label{eq: Q3}
		\begin{split}
		\mathrm{Cyc}_6 \big(  [r_2(x_1&,x_2,x_3,x_4),r_2(x_4,x_5,x_6,x_1)]_{2,1} \\
		 - [r_2&(x_1,x_2,x_4,x_5)]_3 + 2 [r_2(x_1,x_3,x_4,x_5)]_3 \big)  \\
		  &= 4 \, [r_3(x_1,x_2,x_3,x_4,x_5,x_6)]_3 - 6\, [1]_3\ .
		\end{split}
		\end{equation}
		\end{defi}
		
		Les générateurs $[x,y]_{2,1}$ vont jouer le rôle de polylogarithmes multiples (ou corrélateurs) motiviques de profondeur $2$. Comme le suggère l'égalité $B_3(F)=C_3(F)$, ils s'écrivent en fonction des générateurs $[x]_3$ grâce à la spécialisation suivante (au diviseur $\{x_1=x_3\}\simeq \overline{\mathcal{M}}_{0,5}$ du bord de $\overline{\mathcal{M}}_{0,6}$) de la relation~\eqref{eq: Q3}:
		$$[x,y]_{2,1} = [1-x^{-1}]_3  +[1-y^{-1}]_3 + [y/x]_3 +[(1-y)/(1-x)]_3 - [x(1-y)/y(1-x)]_3 - [1]_3\ .$$
		En substituant cette dernière expression dans~\eqref{eq: Q3} on obtient la relation (\og à 22 termes \fg{}) découverte par \textcite{goncharovseattle, goncharovgeometry} et qui est à l'origine de la relation \og à 840 termes \fg{}~\eqref{eq: trirapport}.
		Il est notable que la relation~\eqref{eq: Q3} est liée à la géométrie des configurations de $6$ points dans $\mathbb{P}^1$ alors que~\eqref{eq: trirapport} concerne les configurations de $7$ points dans $\mathbb{P}^2$.
		
		\begin{defi}
		On définit $C_4(F)$ comme le quotient de l'espace vectoriel engendré par des symboles $[x]_4$ et $[x,y]_{3,1}$ pour $x,y\in \mathbb{P}^1(F)$ par les relations suivantes\footnote{\textcite{goncharovrudenko} incluent dans leur définition de $C_4(F)$ les relations entre les générateurs $[x]_4$ qui viennent des relations tétralogarithmiques $\mathcal{R}_4(F)$, tout en conjecturant que cela ne change pas la définition. Cet ajout ne sera pas nécessaire en ce qui nous concerne.}, pour $(x_1,\ldots,x_7)\in \overline{\mathcal{M}}_{0,7}(F)$ :
			\begin{equation}\label{eq: Q4}
			\begin{split}
			\mathrm{Cyc}_7\big( & - [r_2(x_1,x_2,x_3,x_4),r_2(x_4,x_6,x_7,x_1)]_{3,1}  \\
			& + [r_2(x_1,x_2,x_3,x_4),r_2(x_4,x_5,x_7,x_1)]_{3,1} \\
			&-[r_2(x_1,x_2,x_3,x_4),r_2(x_4,x_5,x_6,x_1)]_{3,1} \\
			&+ [r_2(x_1,x_2,x_4,x_6)]_4 + [r_3(x_1,x_2,x_3,x_4,x_5,x_6)]_4\big) = 0\ .\\
			\end{split}
			\end{equation}
		On définit $B_4(F)\subset C_4(F)$ comme le sous-espace engendré par les symboles $[x]_4$.
		\end{defi}
		
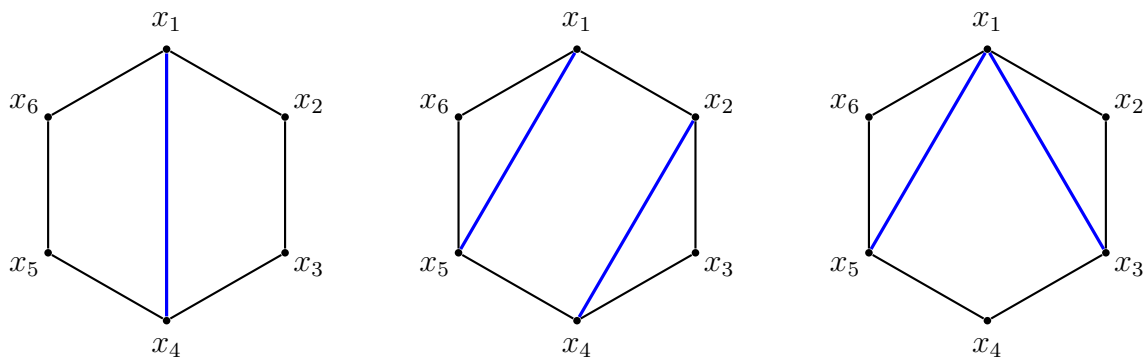
\begin{figure}[h!]
		\begin{center}
		\begin{tikzpicture}

		\def \r {1.8} ;
		
		\foreach \k in {1,...,6} \draw (150-60*\k:\r) node[circle, inner sep = 1pt, fill = black] (\k)  {};
		\foreach \k in {1,...,6} \node at (150-60*\k:1.2*\r) {$x_\k$};
		\draw[thick] (1) -- (2) -- (3) -- (4) -- (5) -- (6) -- (1); 
		\draw[very thick, blue] (1) -- (4);
		
		\begin{scope}[shift = {(3*\r,0)}]
		\foreach \k in {1,...,6} \draw (150-60*\k:\r) node[circle, inner sep = 1pt, fill = black] (\k)  {};
		\foreach \k in {1,...,6} \node at (150-60*\k:1.2*\r) {$x_\k$};
		\draw[thick] (1) -- (2) -- (3) -- (4) -- (5) -- (6) -- (1); 
		\draw[very thick, blue] (1) -- (5);
		\draw[very thick, blue] (2) -- (4);
		\end{scope}
		
		\begin{scope}[shift = {(6*\r,0)}]
		\foreach \k in {1,...,6} \draw (150-60*\k:\r) node[circle, inner sep = 1pt, fill = black] (\k)  {};
		\foreach \k in {1,...,6} \node at (150-60*\k:1.2*\r) {$x_\k$};
		\draw[thick] (1) -- (2) -- (3) -- (4) -- (5) -- (6) -- (1); 
		\draw[very thick, blue] (1) -- (5);
		\draw[very thick, blue] (1) -- (3);
		\end{scope}

		\end{tikzpicture}
		\end{center}
		\caption{Les arguments de $r_2$ dans la relation~\eqref{eq: Q3} sont les sommets des quadrilatères figurant dans les dissections d'un hexagone.}\label{fig: 1}
		\end{figure}

		Les relations~\eqref{eq: Q3} et~\eqref{eq: Q4} ont une origine commune qui sera l'objet du \S\ref{par: cluster}. Contentons-nous pour l'instant de noter que les arguments de $r_2$ et $r_3$ dans ces relations s'interprètent comme des dissections d'un polygone à 6 côtés (respectivement à 7 côtés) comme dans la figure~\ref{fig: 1}.
		
	\subsection{La structure de la cogèbre de Lie motivique en poids $\leq 4$}
	
		Le théorème suivant est prouvé par \textcite{goncharovrudenko}.
	
		\begin{theo}
		\begin{enumerate}[1)]
		\item Les formules suivantes munissent $C_{\leq 4}(F)$ d'une structure de cogèbre de Lie graduée :
		$$\delta[x]_2 = -(1-x)\wedge x \; ,\; \delta[x]_3 = [x]_2\wedge x \; ,\; \delta[x]_4 = [x]_3\wedge x\; , \; $$
		\begin{equation*}
		\begin{split}
		 \delta[x,y]_{2,1} =  [(1-y)/(1-x)]_2 \wedge (y/&x) + [y/x]_2\wedge (1-y)/(1-x) \\
		 & + [x]_2\wedge (1-y^{-1}) + [y]_2\wedge (1-x^{-1})\ ,
		\end{split}
		\end{equation*}
		\begin{equation*}
		\begin{split}
		 \delta[x,y]_{3,1} =  [x,y]_{2,1}\wedge (x/y) & + [x/y]_3\wedge (1-x)/(1-y) 
		 \\ & + [x]_3\wedge (1-y^{-1}) - [y]_3\wedge (1-x^{-1}) - [x]_2\wedge [y]_2\ .
		\end{split}
		\end{equation*}
		\item Les formules suivantes définissent un morphisme de cogèbres de Lie graduées $C_{\leq 4}(F)\to \mathcal{C}_{\leq 4}(F)$
		$$[x]_n\mapsto -\Cor^{\mathcal{C}}(\underbrace{0,\ldots,0}_{n-1},1,x) = \Li_n^{\mathcal{C}}(x)\; \hbox{ pour } n=2,3,4\ ,$$
		$$[x,y]_{2,1} \mapsto \Cor^{\mathcal{C}}(0,x,1,y)\ ,$$
		$$[x,y]_{3,1}\mapsto -\Cor^{\mathcal{C}}(0,0,x,1,y)\ .$$
		\end{enumerate}
		\end{theo}

		C'est la deuxième partie du théorème qui justifie les définitions de $C_{\leq 4}(F)$ et  du cocrochet. Nous expliquerons dans la prochaine section comment Goncharov et Rudenko sont arrivés à ces formules par la combinatoire des dissections de polygones. On renvoie le lecteur aux travaux récents de \textcite{rudenkodepth} et de \textcite{charltonganglradchenkorelations} pour des généralisations en poids supérieur.

		Le théorème suivant, qui est une version \og symbolique \fg{} du théorème~\ref{theo: gangl}, est prouvé par \textcite{goncharovrudenko}.

		\begin{theo}
		On a une suite exacte courte :
		$$0\To B_4(F) \To C_4(F) \stackrel{\delta_{2,2}}{\To} C_2(F)\wedge C_2(F) \To 0 \ .$$
		\end{theo}
		
		On note d'abord que la première flèche est injective par
                définition et que la deuxième flèche est surjective puisque
                $\delta_{2,2}[x,y]_{3,1}=-[x]_2\wedge [y]_2$ pour tous les $x,y\in
                \mathbb{P}^1(F)$. Il reste à montrer l'exactitude au milieu ;
                comme dans la discussion du \S\ref{par: relation 5 termes
                  poids 4}, cela revient à montrer que le morphisme
                $C_2(F)\wedge C_2(F)\To C_4(F)/B_4(F)$
donné                par $[x]_2\wedge [y]_2\mapsto -[x,y]_{3,1}$ est bien défini, c'est-à-dire que le symbole $[x,y]_{3,1}$ vérifie $[x,y]_{3,1}+[y,x]_{3,1}=0$ et la relation à~$5$ termes en la première variable modulo $B_4(F)$. Ces deux relations sont prouvées en combinant des spécialisations de la relation générique~\eqref{eq: Q4} à des strates de $\overline{\mathcal{M}}_{0,7}$ bien choisies.

		On peut montrer que la restriction à $B_{\leq 4}(F)$ du morphisme $C_{\leq 4}(F)\to \mathcal{C}_{\leq 4}(F)$ se factorise par un morphisme $B_{\leq 4}(F)\to \mathcal{B}_{\leq 4}(F)$. 
		On obtient alors le diagramme commutatif de complexes suivant, où $\sim$ désigne un quasi-isomorphisme :
		
		\begin{equation}\label{eq: diag qis crucial}
		\begin{gathered}\xymatrixcolsep{3pc}\xymatrixrowsep{3pc}\xymatrix{
		 \mathcal{B}^\bullet(F,4) \ar[r] & \mathrm{CE}^\bullet(\mathcal{C}(F))_4\  \\ 
		 B^\bullet(F,4) \ar[r]^-{\sim}\ar[u]  & \mathrm{CE}^\bullet(C(F))_4 \  . \ar[u]
		}\end{gathered}
		\end{equation}
		
	\subsection{Compléments sur la preuve}
	
		En plus du quasi-isomorphisme~\eqref{eq: diag qis crucial}, un ingrédient important de la preuve du cas~\mbox{$n=4$} de la conjecture de Zagier par \textcite{goncharovrudenko} est la construction d'un morphisme de complexes
		$$BG^\bullet(F,4) \To \mathrm{CE}^\bullet(C(F))_4\ ,$$
		où $BG^\bullet(F,4)$ est le complexe bigrassmannien introduit au \S\ref{par: strategie premiere}. La composante cruciale est un morphisme $G_7(F,4)\to C_4(F)$, qui associe un élément de $C_4(F)$ à une configuration $(a_0,\ldots,a_7)$ de $8$ points de $F^4$ en position générale. Ce morphisme est une version du polylogarithme grassmannien de \textcite{goncharovarakelov, goncharovsimple} et est construit de manière implicite en se basant sur des travaux antérieurs de \textcite{goncharovgeometrytrilog}.
		
		Peu après, \textcite{charltonganglradchenkograssmannian} ont donné une formule explicite pour ce morphisme, poussé dans $\mathrm{CE}^\bullet(\mathcal{C}(F))_4$, en termes de $\mathrm{I}_{3,1}^{\mathcal{C}}$ et $\Li_4^{\mathcal{C}}$. Il devrait donc être possible, quitte à vérifier des compatibilités avec les relations définissant $C(F)$, de se servir de cette formule pour rendre la preuve de \textcite{goncharovrudenko} effective.
		
		De plus, \textcite{charltonganglradchenkograssmannian} écrivent explicitement un représentant dans $\mathcal{B}_4(F)$ du même cocycle pour en déduire une équation fonctionnelle du tétralogarithme, dans l'esprit de~\eqref{eq: five term birapport} et~\eqref{eq: trirapport}. Elle s'écrit dans $\mathcal{B}_4(F)$ sous la forme prévue
		$$\sum_{i=0}^8 (-1)^i R_4(a_0,\ldots,\widehat{a_i},\ldots,a_8) =0\ ,$$
		où $R_4(a_1,\ldots,a_8)\in \QQ[F^\times]$ est un invariant explicite mais compliqué d'une configuration de $8$ points dans $\mathbb{P}^3(F)$, dont l'interprétation géométrique n'est pas encore claire.

\section{Le rôle des structures amassées et de la combinatoire des dissections}\label{par: cluster}

	Les structures amassées, introduites par \textcite{fominzelevinsky1}, sont des objets de combinatoire algébrique (et géométrique) qui interviennent dans de nombreux domaines. \textcite{fockgoncharov} ont mis en lumière l'importance de ces structures dans l'étude du dilogarithme et de ses aspects motiviques, et \textcite{goncharovrudenko} ont incorporé les polylogarithmes (multiples) supérieurs dans ce contexte. Ce sont ces liens, et notamment la combinatoire des dissections de polygones, qui expliquent la forme des relations~\eqref{eq: Q3} et~\eqref{eq: Q4} dont il a été question plus haut.

	\subsection{Structures amassées et dilogarithme}
	
		Fixons un entier $n$. Une \emph{graine} $(B,u)$ est la donnée d'une matrice antisymétrique à coefficients entiers $B=(b_{ij})_{1\leq i,j\leq n}$ et d'une famille $u=(u_1,\ldots,u_n)$ d'éléments du corps $\QQ(x_1,\ldots,x_n)$ qui l'engendrent librement. On peut la représenter par un carquois (graphe orienté) dont les sommets sont les entiers de $1$ à $n$, décorés par les $u_i$, et où $b_{ij}$ arêtes vont de $i$ vers $j$ si $b_{ij}\geq 0$. On appelle $B$ la \emph{matrice d'échange}, $u$ l'\emph{amas} et les $u_i$ les \emph{variables d'amas}.
		
		Pour un sommet $k\in \{1,\ldots,n\}$ la \emph{mutation} de la graine $(B,u)$ suivant le sommet $k$ est la graine $(B',u')$ où
		$$b'_{ij} = 
		\begin{cases}
		-b_{ij} & \hbox{ si } i=k \hbox{ ou } j=k \\ 
		b_{ij} +|b_{ik}|b_{kj} & \hbox{ si } b_{ik}b_{kj}>0 \\
		b_{ij} & \hbox{ sinon}
		\end{cases}
		$$
		et les variables $u'_i$ sont données par $u'_i=u_i$ pour $i\neq k$ et par la relation d'échange
		$$u_ku'_k =\prod_{i\,|\, b_{ik}>0}u_	i^{b_{ik}} + \prod_{j\,|\, b_{kj}>0} u_j^{b_{kj}}\ .$$
		
		\begin{exem}
		Voici un exemple de mutation (suivant le sommet de gauche) :
		$$x_1\;\stackrel{}{\bullet}\To\stackrel{}{\bullet}\;x_2 \hspace{1cm} \leadsto \hspace{1cm} \textstyle\frac{1+x_2}{x_1} \;\stackrel{}{\bullet} \longleftarrow \stackrel{}{\bullet}\; x_2$$
		Si l'on continue à appliquer des mutations suivant le sommet de droite, puis de gauche, puis de droite, puis de gauche, on retombe (après $5$ mutations donc) sur la graine de départ, à renumérotation des sommets près : $x_2 \;\stackrel{}{\bullet}\longleftarrow \stackrel{}{\bullet}\; x_1$.
		\end{exem}
		
		En partant d'une graine et en appliquant successivement des mutations suivant tous les sommets possibles, on produit ce qu'on appellera improprement une \emph{structure amassée}. Il arrive, comme dans l'exemple précédent, qu'on ne produise ainsi qu'un nombre fini de graines différentes. Ce phénomène de finitude est rare et correspond, d'après \textcite{fominzelevinsky2}, aux graines dont une mutation a un carquois sous-jacent qui est une orientation d'un carquois de Dynkin (type ADE). Dans ces cas, la structure combinatoire des mutations peut être modélisée par un graphe appelé \emph{graphe d'échange}, où les sommets sont les graines (à renumérotation des sommets près) et les arêtes indiquent les mutations. D'après \textcite{chapotonfominzelevinsky}, le graphe d'échange est le $1$-squelette d'un polytope convexe appelé associaèdre généralisé. En type $A$ il s'agit de l'associaèdre classique, découvert par Stasheff. L'exemple que nous venons de traiter correspond au type $A_2$, et son graphe d'échange est un cycle de longueur $5$, qui est le bord d'un pentagone.
	
		L'importance des structures amassées dans l'étude du dilogarithme et de ses aspects motiviques a été mise en lumière par \textcite{fockgoncharov}. Pour une graine~$(B,u)$ donnée on définit :
		$$W = \frac{1}{2}\, \sum_{i,j} b_{ij} \, u_i\wedge u_j \;\; \in \Lambda^2(\QQ(x_1,\ldots,x_n)^\times_\QQ)\ .$$
		Si $(B',u')$ est obtenue par mutation du sommet $k$, on a la relation fondamentale
		\begin{equation}\label{eq: def X}
		W' - W = - (1+X)\wedge X \;\;\;\hbox{ avec }\;\; X = \prod_{j=1}^n u_j^{b_{kj}}\ .
		\end{equation}
		On l'écrit sous la forme
                plus suggestive
		\begin{equation}\label{eq: mutation differential}
		W'-W = \delta_2[-X]_2\ ,
		\end{equation}
		où $[-X]_2$ est vu dans le groupe de Bloch $\mathcal{B}_2(\QQ(x_1,\ldots,x_n))$. Cela implique qu'on a l'égalité $[W]=[W']$ dans $K_2(\QQ(x_1,\ldots,x_n))$ : une classe canonique dans le $K_2$ est attachée à la structure amassée.
		
		En revenant au cas du type $A_2$, on voit qu'on peut exprimer la discussion précédente sous la forme du diagramme commutatif suivant :
		\begin{equation*}\label{eq: morphisme complexes A2}
		\begin{gathered}
		\xymatrixrowsep{1.5cm}\xymatrixcolsep{3.5pc}\xymatrix{
		0 \ar[r] & C_2\ar[d]\ar[r]^-{\partial} & C_1 \ar[d]\ar[r]^-{\partial} & \ar[d]C_0 \ar[r] & 0\ \\
		0 \ar[r] & 0 \ar[r] & \mathcal{B}_2(\QQ(x_1,x_2)) \ar[r]^-{\delta_2} & \Lambda^2(\QQ(x_1,x_2)^\times)_\QQ \ar[r] & 0\ .
		}
		\end{gathered}
		\end{equation*}
		La première ligne est le complexe cellulaire du pentagone. La flèche verticale de droite envoie un sommet du pentagone, c'est-à-dire une graine, vers l'élément $W$ correspondant. La flèche verticale du milieu envoie une arête du pentagone, c'est-à-dire une mutation, vers l'élément $[-X]_2$ correspondant. La commutativité du carré de droite est équivalente à~\eqref{eq: mutation differential}. La commutativité du carré de gauche est équivalente à l'égalité
		$$[-x_1]_2 + \left[-\frac{1+x_2}{x_1}\right]_2 + \left[-\frac{1+x_1+x_2}{x_1x_2}\right]_2 + \left[-\frac{1+x_1}{x_2}\right]_2 + [-x_2]_2 = 0$$
		dans $\mathcal{B}_2(\QQ(x_1,x_2))$, ce qui est une réécriture\footnote{Il suffit pour s'en convaincre de poser $x_1=-x$, $x_2=-1+y$, et d'utiliser les relations $[1-u]_2=-[u]_2$ et $[u^{-1}]_2=-[u]_2$. Notons qu'on peut justifier la commutativité du carré de gauche, et donc redécouvrir la relation à $5$ termes, par un argument abstrait. En effet, $\delta_2$ est injective puisque d'après la remarque \ref{rema: suslin bloch} son noyau est un quotient de $\mathrm{gr}^2_\gamma K_3(\QQ(x_1,x_2))_\QQ$, qui par localisation est isomorphe à~$\mathrm{gr}^2_\gamma K_3(\QQ)_\QQ$, et donc nul par le théorème \ref{theo: borel}. La commutativité du carré de droite et le fait que $\partial\circ\partial =0$ impliquent donc la commutativité du carré de gauche.
		} de la relation à $5$ termes du dilogarithme~\eqref{eq: five term abstract}.
		
		Une idée importante de \textcite{goncharovrudenko} est que les liens entre structures amassées et dilogarithmes devraient s'étendre en poids supérieur. On se contentera de développer ces idées dans le cadre du type $A$, alors que la vision de Goncharov--Rudenko est plus générale. Il est d'ailleurs connu que les autres types jouent un rôle dans les équations fonctionnelles des polylogarithmes (on trouve chez \textcite{ggsvv} une équation fonctionnelle pour le trilogarithme dans le cadre de la structure amassée de type $D_4$).

	\subsection{Structures amassées en type $A$}
	
 		Rappelons que l'associaèdre $\mathcal{A}_n$ est un polytope convexe de dimension $n-3$ dont le treillis des faces est en bijection avec le treillis des dissections d'un polygone $\Pi_n$ à $n$ côtés, la codimension d'une face étant le nombre de diagonales dans la dissection correspondante. Les sommets de $\mathcal{A}_n$ sont donc en bijection avec les triangulations de $\Pi_n$. La famille des associaèdres a une structure opéradique : à une dissection de $\Pi_n$ en des polygones $\Pi_{n_1},\ldots,\Pi_{n_r}$ correspond une face de $\mathcal{A}_n$ isomorphe au produit $\mathcal{A}_{n_1}\times\cdots \times \mathcal{A}_{n_r}$. 
 		
 		\begin{figure}[h!]
		\begin{center}\includegraphics[width=0.7\textwidth]{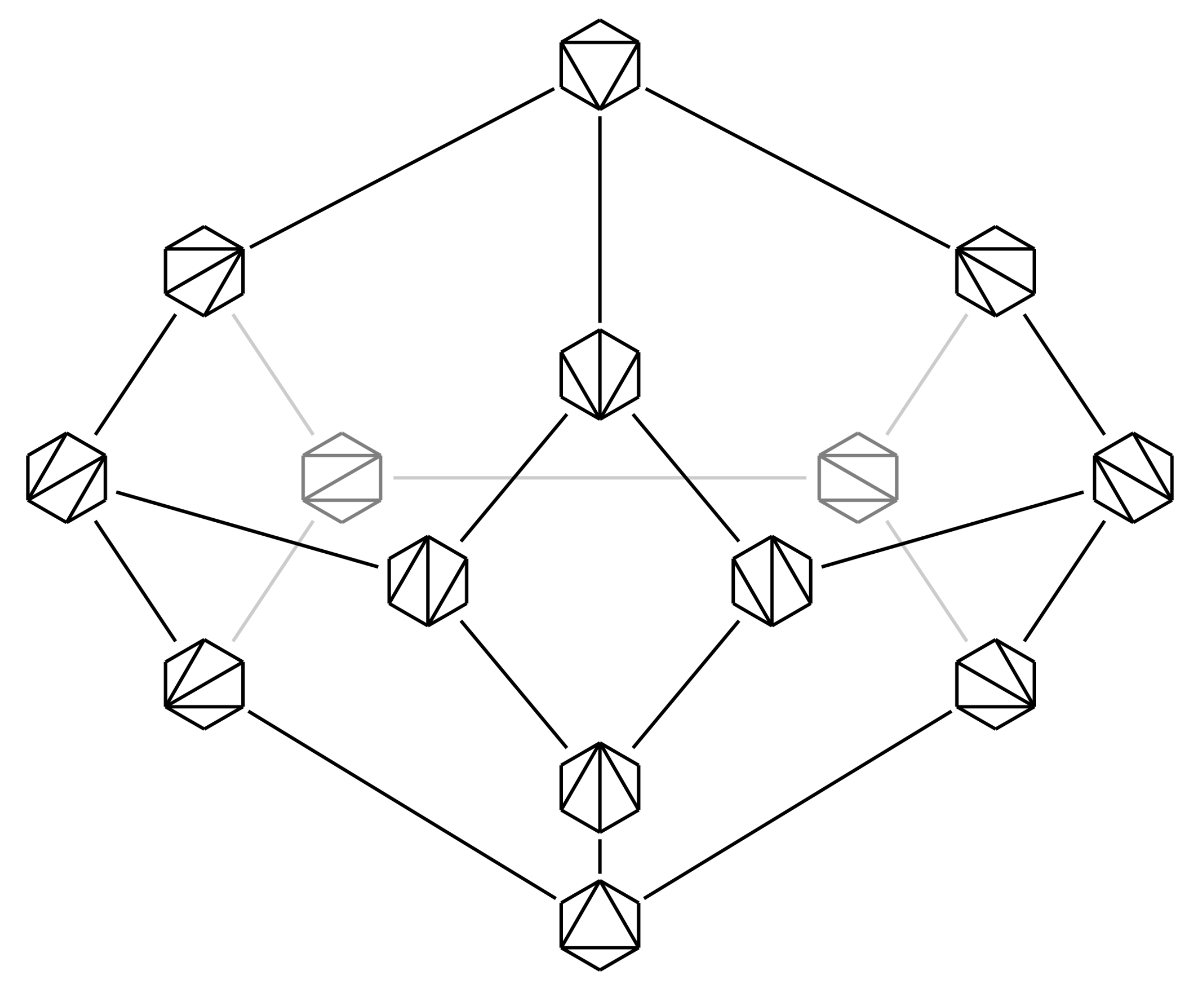}\end{center}
		\caption{L'associaèdre $\mathcal{A}_6$ a $3$ faces carrées et $6$ faces pentagonales (figure empruntée à \textcite{fominzelevinskyCDM}).}\label{fig: associaedre}
		\end{figure}
 		
 		Le $1$-squelette de $\mathcal{A}_n$ est le graphe d'échange d'une structure amassée (de type $A_{n-3}$) qu'on décrit maintenant. Notons $|v,w|$ le déterminant de deux vecteurs $v,w$ du plan. Soit $Y_n$ la variété sur $\QQ$ définie comme le quotient par $\mathrm{SL}_{2}$ de l'espace des $n$-uplets $(v_1,\ldots,v_n)$ de points du plan en position générale (c'est-à-dire tels que $|v_i,v_j|\neq 0$ si~$i\neq j$). Via les déterminants $|v_i,v_j|$, on peut voir $Y_n$ comme la sous-variété du tore $\mathbb{G}_{m}^{n(n-1)/2}$ définie par les équations de Plücker. On peut aussi voir $Y_n$ comme l'ouvert du cône affine sur la Grassmannienne $\mathrm{Gr}(2,n)$ formé des $2$-plans dans un espace vectoriel de dimension $n$ qui sont en position générale par rapport à $n$ hyperplans de coordonnées.  Notons $F_n$ le corps des fonctions de $Y_n$, qui est non canoniquement isomorphe à $\QQ(x_1,\ldots,x_{2n-3})$. \`{A} une dissection de $\Pi_n$ en des polygones $\Pi_{n_1},\ldots,\Pi_{n_r}$ est associé un morphisme de $Y_n$ vers $Y_{n_1}\times\cdots \times Y_{n_r}$ dont les composantes sont des morphismes d'oubli, et on a donc un morphisme $F_{n_1}\otimes \cdots \otimes F_{n_r}\To F_n$.
 		
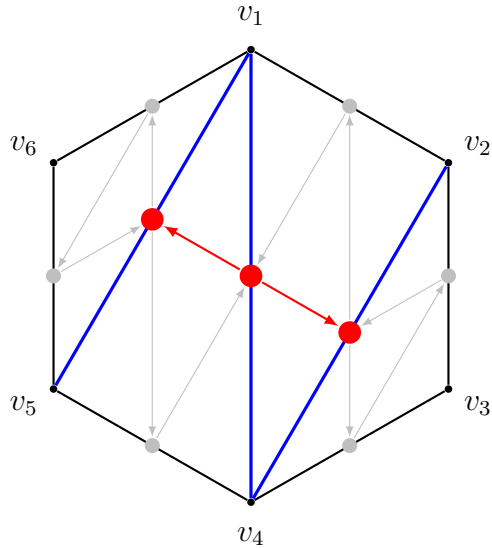
\begin{figure}[h!]
		\begin{center}
		\begin{tikzpicture}

		\def \r {3} ;
		\foreach \k in {1,...,6} \draw (150-60*\k:\r) node[circle, inner sep = 1pt, fill = black] (\k)  {};
		\foreach \k in {1,...,6} \node at (150-60*\k:1.15*\r) {$v_\k$};
		\draw[thick] (1) -- (2) node[circle,inner sep = 2pt, fill=lightgray, midway] (mid12) {}; 
		\draw[thick] (2) -- (3) node[circle,inner sep = 2pt, fill=lightgray, midway] (mid23) {}; 
		\draw[thick] (3) -- (4) node[circle,inner sep = 2pt, fill=lightgray, midway] (mid34) {}; 
		\draw[thick] (4) -- (5) node[circle,inner sep = 2pt, fill=lightgray, midway] (mid45) {};
		\draw[thick] (5) -- (6) node[circle,inner sep = 2pt, fill=lightgray, midway] (mid56) {}; 
		\draw[thick] (6) -- (1) node[circle,inner sep = 2pt, fill=lightgray, midway] (mid61) {};
		
		\draw[very thick, blue] (1) -- (4) node[circle, inner sep = 3pt, fill=red, midway] (mid14) {} ;
		\draw[very thick, blue] (1) -- (5) node[circle, inner sep = 3pt, fill=red, midway] (mid15) {} ;
		\draw[very thick, blue] (2) -- (4) node[circle, inner sep = 3pt, fill=red, midway] (mid24) {} ;
		
		\draw[->,>=latex, red, thick] (mid14) -- (mid24) ;
		\draw[->,>=latex, red, thick] (mid14) -- (mid15) ;
		\draw[->,>=latex, lightgray] (mid24) -- (mid12) ;
		\draw[->,>=latex, lightgray] (mid12) -- (mid14) ;
		\draw[->,>=latex, lightgray] (mid24) -- (mid34) ;
		\draw[->,>=latex, lightgray] (mid34) -- (mid23) ;
		\draw[->,>=latex, lightgray] (mid23) -- (mid24) ;
		\draw[->,>=latex, lightgray] (mid15) -- (mid45) ;
		\draw[->,>=latex, lightgray] (mid45) -- (mid14) ;
		\draw[->,>=latex, lightgray] (mid15) -- (mid61) ;
		\draw[->,>=latex, lightgray] (mid61) -- (mid56) ;
		\draw[->,>=latex, lightgray] (mid56) -- (mid15) ;
		\end{tikzpicture}
		\end{center}
		\caption{Le carquois associé à une triangulation.}\label{fig: triangulation carquois}
		\end{figure}

		On oriente le polygone $\Pi_n$ et on étiquette ses sommets avec des symboles $v_1,\ldots,v_n$ de manière cyclique. On associe une graine à une triangulation de $\Pi_n$ de la manière suivante, illustrée par la figure~\ref{fig: triangulation carquois}. Définissons un carquois en plaçant un sommet sur chacune des $n-3$ diagonales de la triangulation et sur chacun des $n$ côtés de $\Pi_n$ (ces $n$ derniers sommets sont considérés \emph{gelés}, ce qui signifie qu'ils ne donneront pas lieu à des mutations). On décore par la variable d'amas $|v_i,v_j|\in F_n$ un sommet qui est sur un segment reliant $v_i$ et $v_j$, avec $i<j$. Chaque triangle de la triangulation donne lieu à trois sommets du carquois, qu'on relie par trois arêtes orientées dans le sens direct.

\begin{figure}[h]
		\begin{center}
		
		\begin{tikzpicture}
		
		\def \a {4} ;
		\draw[thick] (0,0) rectangle (\a,\a) ;
		\draw[very thick, blue] (\a,0) -- (0,\a) node[circle, inner sep = 3pt, fill=red, midway] {};
		\filldraw[black] (0,\a) circle (\a/4 pt) node[anchor=south east] {$v_1$} ;
		\filldraw[black] (\a,\a) circle (\a/4 pt) node[anchor=south west] {$v_2$} ;
		\filldraw[black] (\a,0) circle (\a/4 pt) node[anchor=north west] {$v_3$} ;
		\filldraw[black] (0,0) circle (\a/4 pt) node[anchor=north east] {$v_4$} ;
		
		\begin{scope}[shift={(1.5*\a,0.5*\a)}]
		\draw node {$\leadsto$} ;
		\end{scope}
		
		\begin{scope}[shift = {(2*\a, 0)}]
		\draw[thick] (0,0) rectangle (\a,\a) ;
		\draw[very thick, blue] (0,0) -- (\a,\a) node[circle, inner sep = 3pt, fill=red, midway] {} ;
		\filldraw[black] (0,\a) circle (\a/4 pt) node[anchor=south east] {$v_1$} ;
		\filldraw[black] (\a,\a) circle (\a/4 pt) node[anchor=south west] {$v_2$} ;
		\filldraw[black] (\a,0) circle (\a/4 pt) node[anchor=north west] {$v_3$} ;
		\filldraw[black] (0,0) circle (\a/4 pt) node[anchor=north east] {$v_4$} ;
		\end{scope}
		
		\end{tikzpicture}
		\end{center}
		\caption{Un \og flip \fg{} d'une diagonale dans un quadrilatère.}\label{fig: triangulation flip}
		\end{figure}
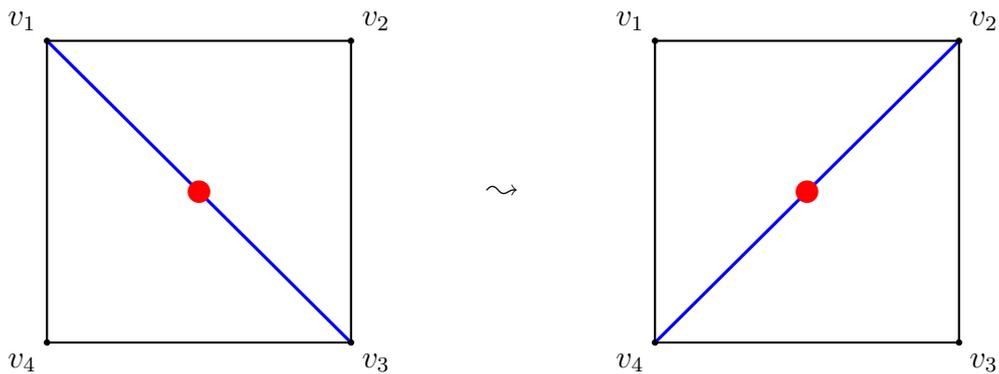

		Soient $T$ et $T'$ deux triangulations qui sont reliées par un \og flip \fg{} d'une diagonale dans un quadrilatère, comme dans la figure~\ref{fig: triangulation flip}. On voit grâce aux relations de Plücker que les graines associées à $T$ et $T'$ sont des mutations l'une de l'autre suivant le sommet correspondant au quadrilatère.  On engendre, en itérant les \og flips \fg{} (mutations suivant les sommets non gelés), une structure amassée avec un nombre fini de graines dont le graphe d'échange est le $1$-squelette de l'associaèdre $\mathcal{A}_n$.

		Par la construction du paragraphe précédent, on associe un élément $W_T\in \Lambda^2(F_n^\times)_\QQ$ à toute triangulation $T$ de $\Pi_n$. Concrètement, on associe au triangle ($n=3$) l'élément 
		$$W_\triangle = |v_1,v_2|\wedge |v_2,v_3| + |v_2,v_3|\wedge |v_1,v_3| + |v_1,v_3|\wedge |v_1,v_2|  \;\; \in \Lambda^2(F_3^\times)_\QQ\ ,$$
		et à une triangulation $T$ la somme des contributions de chaque triangle :
		\begin{equation}\label{eq: W triangles}
		W_T = \sum_{\triangle\in T}W_\triangle\ .
		\end{equation}
		
	\subsection{Polylogarithmes amassés}
	
		\subsubsection{La conjecture}
	
		\textcite{goncharovrudenko} partent de l'idée, simple mais féconde, d'écrire la relation~\eqref{eq: W triangles} sous une forme multiplicative, et donc plus naturellement compatible à la structure opéradique,
		\begin{equation}\label{eq: exp W triangles}
		\exp(W_T) = \prod_{\triangle\in T} \exp(W_\triangle)\ ,
		\end{equation}
		dans l'algèbre extérieure $\Lambda(F_n^\times)_\QQ$ (il n'y a qu'un nombre fini de termes puisque $\exp(W_\triangle)=1+W_\triangle$). C'est ce passage à l'exponentielle qui permet de passer du dilogarithme aux poids supérieurs. Afin de formaliser cela, il sera commode de faire intervenir la cogèbre de Lie motivique du corps $F_n$, dont l'existence est malheureusement conditionnelle à la conjecture de Beilinson--Soulé pour $F_n$, qui n'est pas prouvée. On pourrait écrire des énoncés inconditionnels (mais moins propres) au prix de remplacer les $F_n$ par un corps de nombres $F$ et les variables $v_i$ par des vecteurs suffisamment génériques de $F^2$.
		
		Passons en notation homologique et considérons les complexes
		$$\mathrm{CE}_\bullet(\mathcal{C}(F_n))_k = \mathrm{CE}^{k-\bullet}(\mathcal{C}(F_n))_k \;\;\hbox{ et }\;\; \mathrm{CE}_\bullet(\mathcal{C}(F_n)) = \bigoplus_{k\geq 0}\mathrm{CE}_\bullet(\mathcal{C}(F_n))_k\ .$$
		On note $C_\bullet(\mathcal{A}_n)$ le complexe cellulaire de l'associaèdre $\mathcal{A}_n$, dont une base est donnée par les dissections de $\Pi_n$. La conjecture suivante est due à \textcite{goncharovrudenko}.
		
		\begin{conj}\label{conj: cluster polylog}
		Il existe des morphismes de complexes
		$$\alpha_n : C_\bullet(\mathcal{A}_n) \To \mathrm{CE}_\bullet(\mathcal{C}(F_n))$$ pour $n\geq 3$, appelés \emph{polylogarithmes amassés}, qui vérifient :
		\begin{enumerate}[---]
		\item Initialisation : le morphisme  $\alpha_{3} : C_0(\mathcal{A}_3) \to \mathrm{CE}_0(\mathcal{C}(F_3)) =\Lambda(F_3^\times)_\QQ$ est donné par 
		$$\alpha_3(\triangle) = \exp(W_\triangle)\ .$$
		\item Compatibilité aux structures opéradiques : pour toute dissection $D_1\times D_2$ de $\Pi_n$ obtenue en recollant une dissection $D_1$ de $\Pi_{n_1}$ et une dissection $D_2$ de $\Pi_{n_2}$ avec $n_1+n_2=n+1$, on a l'égalité :
		\begin{equation*}\label{eq: compatibilite polylog cluster}
		\alpha_n(D_1\times D_2) = \alpha_{n_1}(D_1)\,\alpha_{n_2}(D_2)\ ,
		\end{equation*}
		où le membre de droite s'interprète via les inclusions naturelles de $F_{n_1}$ et $F_{n_2}$ dans~$F_n$.
		\end{enumerate}
		\end{conj}
		
		Grâce à la structure opéradique, les images par $\alpha_n$ des générateurs de $C_{n-3}(\mathcal{A}_n)$, pour tout $n\geq 3$, déterminent uniquement la collection des morphismes $\alpha_n$. L'idée est de découvrir les $\alpha_n$ par une procédure récursive.
		
		\begin{rema}
		D'un point de vue opéradique, la conjecture~\ref{conj: cluster polylog} prévoit l'existence d'un diagramme commutatif d'opérades (différentielles graduées, cycliques non symétriques) :
		$$\xymatrixcolsep{4pc}\xymatrixrowsep{4pc}\xymatrix{
		\mathrm{As}_\infty \ar@{-->}[d]_-{\alpha} \ar[r]^-{\sim} & \mathrm{As}\ \ar[d]^{\exp(W)} \\
		\mathrm{CE}(\mathcal{C}) \ar[r] & K^{\mathrm{M}}	\ .	
		}$$
		Ce diagramme fait intervenir l'opérade $\mathrm{As}$ des algèbres associatives et sa résolution canonique $\mathrm{As}_\infty$, l'opérade des algèbres associatives à homotopie près, données en arité~$n$ par $\mathrm{As}(n)=H_0(\mathcal{A}_n)=\QQ$ et $\mathrm{As}_\infty(n)=C_\bullet(\mathcal{A}_n)$, où la structure opéradique est induite par les inclusions des faces des associaèdres (voir par exemple \textcite[chapitre 9]{lodayvallette}). On a noté $\mathrm{CE}(\mathcal{C})$ et $K^{\mathrm{M}}$ les opérades données en arité $n$ par $\mathrm{CE}(\mathcal{C})(n)=\mathrm{CE}_\bullet(\mathcal{C}(F_n))$ et $K^{\mathrm{M}}(n)=H_0(\mathrm{CE}_\bullet(\mathcal{C}(F_n)))=K^{\mathrm{M}}(F_n)$ respectivement, avec les structures opéradiques induites par le produit de l'algèbre extérieure. Il ne semble pas y avoir d'argument abstrait d'algèbre homologique qui implique l'existence du morphisme $\alpha$.
		\end{rema}

		\subsubsection{Calculs en basse dimension}
		
		En degré $0$, la compatibilité aux structures opéradiques et~\eqref{eq: exp W triangles} forcent $\alpha_n$ à être définie par $\alpha_n(T) = \exp(W_T)$. En degré $1$, l'égalité~\eqref{eq: mutation differential} motive la définition, pour $D$ une dissection qui contient des triangles et un quadrilatère :
		$$ \alpha_n(D) = - \Li_2^{\mathcal{C}}(-X_D)\otimes \frac{1}{2}(\exp(W_T)+\exp(W_{T'}))\ , $$
		où $T$ et $T'$ sont les deux triangulations de $\Pi_n$ compatibles avec $D$, et $X_D$ est défini comme dans~\eqref{eq: def X}. On vérifie que cette formule est compatible aux structures opéradiques et qu'on a $\delta\alpha_n(D)=\alpha_n(T')-\alpha_n(T) = \alpha_n(\partial D)$. Notons que $-X_D$ est un birapport ; par exemple dans le cas du carré on a en poids $2$ :
		$$\alpha_4(\square)_2 = -\Li_2^{\mathcal{C}}(r_2(x_1,x_2,x_3,x_4))\ ,$$ 
		où l'on note $x_i\in \mathbb{P}^1$ la classe de $v_i$.	
		
		Goncharov et Rudenko étendent la définition de $\alpha_n$ en degré $2$. Pour cela il suffit de traiter le cas d'une dissection $D$ de $\Pi_n$ qui contient des triangles et un pentagone. Notons $D_1,\ldots,D_5$ les dissections obtenues en ajoutant une diagonale à $D$, de telle sorte que la superposition de $D_{i-1}$ et $D_{i}$ (indices modulo $5$) crée une triangulation $T_i$. On remarque que l'élément
		$$\widetilde{W} = W_{D_i} - X_{D_{i-1}}\wedge X_{D_i} \;\; \in \Lambda^2(F_n^\times)_\QQ$$
		est indépendant de $i\in \mathbb{Z}/5\mathbb{Z}$. La formule
		$$\alpha_n(D) = \frac{1}{2}\sum_{i\in\mathbb{Z}/5\mathbb{Z}}\Li_3^{\mathcal{C}}(-X_{D_i})\otimes (X_{D_{i-1}}/X_{D_{i+1}}) \wedge \exp(\widetilde{W})$$
		est alors compatible aux structures opéradiques et vérifie $\delta\alpha_n(D)= \alpha_n(\partial D)$. Par exemple dans le cas du pentagone on a en poids $4$ :
		$$\alpha_5(\pentagon)_4 = \frac{1}{2} \mathrm{Cyc}_5\big(\Li_3^{\mathcal{C}}(r_2(x_1,x_2,x_3,x_4))\otimes r_2(x_1,x_3,x_4,x_2)\big)\ .$$
		
		En degré $3$, Goncharov et Rudenko utilisent la combinatoire des corrélateurs pour découvrir une formule pour $\alpha_6$. Ils prouvent qu'en posant
		\begin{equation*}
		\begin{split}
		\alpha_6(\varhexagon)_4 = \mathrm{Cyc}_6^-\big(-\Cor^{\mathcal{C}}&(0,0,r_2(x_1,x_2,x_3,x_4),1,r_2(x_4,x_5,x_6,x_1)) \\
		& +\Li_4^{\mathcal{C}}(r_2(x_1,x_3,x_4,x_5)) - \Li_4^{\mathcal{C}}(r_3(x_1,x_2,x_3,x_4,x_5,x_6)) \big)\ ,
		\end{split}
		\end{equation*}
		on définit un morphisme de complexes $C_\bullet(\mathcal{A}_6) \to \mathrm{CE}_\bullet(\mathcal{C}(F_6))_4$. 
		
		La structure opéradique permet alors d'induire un morphisme de complexes $C_{\bullet\leq 3}(\mathcal{A}_7)\to \mathrm{CE}_\bullet(\mathcal{C}(F_7))_4$ puisque les faces de dimension $\leq 3$ de l'associaèdre $\mathcal{A}_7$ sont toutes des produits d'associaèdres de dimension inférieure. On a donc l'égalité, dans $\mathrm{CE}_2(\mathcal{C}(F_7))_4 = \mathcal{C}_3(F_7)\wedge (F_7^\times)_\QQ\oplus \mathcal{C}_2(F_7)\wedge \mathcal{C}_2(F_7)$ :
		$$\delta(\alpha_7( \partial( \heptagon ) )_4 ) = \alpha_7(\partial\partial (\heptagon))_4 = 0\ .$$
		Par un argument de rigidité comme dans la remarque~\ref{rema: relations periodes motiviques}, on a donc, pour $F$ un corps de nombres et $x_1,\ldots,x_7\in\mathbb{P}^1(F)$, la relation
		$$\alpha_7(\partial(\heptagon))_4=0$$
		dans $\mathcal{C}_4(F)$. C'est exactement la relation~\eqref{eq: Q4} qui apparaît dans la définition de $C_4(F)$. 
	
	\subsection{Polylogarithmes quadrangulaires}
	
		Concluons ce texte avec un mot sur les travaux de \textcite{rudenkodepth} et la preuve du théorème~\ref{theo: rudenko depth}, qui est une étape vers la conjecture de profondeur de Goncharov. Les techniques utilisent des objets appelés \emph{polylogarithmes quadrangulaires}. Ils sont notés
		$$\mathrm{QLi}_{n,k}(x_0,\ldots,x_{2n+1}) \in \mathcal{C}_{n+k}(F)\ ,$$
		avec $x_0,\ldots,x_{2n+1} \in \mathbb{P}^1(F)$, et sont définis comme des sommes alternées de certains corrélateurs $\mathrm{Cor}^{\mathcal{C}}(x_{i_0},\ldots,x_{i_{n+k}})$. Rudenko montre les deux liens suivants entre polylogarithmes quadrangulaires et polylogarithmes multiples.
		\begin{enumerate}[(a)]
		\item Tout corrélateur motivique de poids pair $2n$ s'écrit comme combinaison linéaire de $\mathrm{QLi}_{n,n}$. Tout corrélateur motivique de poids impair $2n+1$ s'écrit comme combinaison linéaire de $\mathrm{QLi}_{n,n+1}$.
		\item Le polylogarithme quadrangulaire $\mathrm{QLi}_{n,k}$ s'écrit comme combinaison linéaire de polylogarithmes multiples motiviques de profondeur $\leq n$.
		\end{enumerate}
		La partie la plus subtile est la preuve de (b), qui s'appuie sur la combinatoire des quadrangulations (dissection en des quadrilatères) des polygones. La relation exacte entre polylogarithmes quadrangulaires et polylogarithmes amassés reste à préciser.

\printbibliography	
	
\end{document}